\documentclass[11pt]{amsart}
\usepackage{amsfonts,amssymb,stmaryrd,amscd,amsmath,latexsym,amsbsy}

\theoremstyle{plain}
\newtheorem{thm}{Theorem}
\newtheorem*{thrm}{Theorem}

\newtheorem{pro}{Proposition}
\newtheorem{lem}{Lemma}
\newtheorem{cor}{Corollary}

\newtheorem*{sol}{Solution}
\newtheorem*{LRrule}{Littlewood-Richardson rule}

\theoremstyle{definition}
\newtheorem*{defi}{Definition}

\theoremstyle{remark}
\newtheorem{rema}{Remark}

\newcommand{\solu}[1]{\begin{sol}{\bf (\ref{#1})}}
\newcommand{\dd}{\partial}
\newcommand{\gotg}{\mathfrak g}
\newcommand{\gotk}{\mathfrak k}
\title{Jack polynomials for the $BC_n$ root system and
generalized spherical functions}
\author{A.~Oblomkov}
\address{Department of Mathematics, MIT, 77, Massachusetts Ave., Cambridge,
MA 02139, USA.}

\email{oblomkov@math.mit.edu}

\begin{document}

\maketitle
\section*{Introduction}
Functions on a homogeneous space   $G/K$ invariant with respect to
the left action of $K$ are called spherical functions (or
sometimes $K$-spherical). One can also study functions on $G/K$
with values in a representation $V$ of $G$ which are equivariant
with respect to the left action of $K$. This more general class of
functions may be called vector valued spherical functions. The
theory of such functions was developed by Harish-Chandra, Helgason
and other authors \cite{HC},\cite{He},\cite{W}.

In the case when $G=K\times K$ and $K\subset G$ is the diagonal
subgroup, the study of
 vector valued spherical functions is equivalent to the study of
 functions on $K$,
equivariant with respect to   conjugation. When $K$ is reductive
over $\mathbb C$, the Peter-Weyl theorem gives  a description of
the space of conjugacy equivariant functions  as the space spanned
by the vector valued characters of $K$. The article \cite{EK4}
deals with such kind of spherical functions in the case when
$K=SL(n,\mathbb C)$. Namely, the Laplace operator on $K$
restricted to the space of spherical functions can be written in
terms of coordinates along  the maximal torus. If the vector
valued  equivariant functions take values in $V=S^{\kappa
n}\mathbb C^n$, then the resulting operator (the so called radial
part of the Laplacian)
 coincides, up to an obvious conjugation, with
the Sutherland differential operator, which is the Hamiltonian of the
 Calogero-Moser quantum
mechanical system  for the root system $A_{n-1}$ \cite{C},
\cite{OP}, \cite{S}. Differential operators on $K$ corresponding
to the higher Casimir operators can also be written in terms of
coordinates along the maximal torus. These operators are quantum
integrals of the
 Calogero-Moser system. Moreover,
Weyl group invariant eigenfunctions of the
 Sutherland operator
can be expressed as vector valued traces of some intertwining
operators between some particular representations of $K$. By the
results from \cite{HO} these eigenfunctions are essentially Jack
polynomials for the root system $A_{n-1}$ (up to a
Weyl-determinantlike factor).

The main result of this paper   is a representation theoretic
interpretation (in the spirit of the work \cite{EK4}) of the three
parameter family of $BC_n$ Jack polynomials.
 More precisely, we consider
the case of the pair
 $G=GL(m+n,\mathbb C)$ ($m\ge n$),
$K=GL(m,\mathbb C)\times GL(n,\mathbb C)$, slightly modify the
definition of spherical functions and as result we get a similar
to \cite{EK4} theory for the root system $BC_n$. Namely, in this
case the Laplace operator on $G$, written in terms of coordinates
of some torus inside $G$, yields the Sutherland operator for the
root system $BC_n$, and the higher Casimir operators give  quantum
integrals of the corresponding quantum mechanical system.

Furthermore, the restriction to the torus of some special matrix
elements of irreducible finite dimensional representations
$L_\lambda$ of $G$, where $\lambda$ ranges over the set isomorphic
to the cone of dominant integral weights of the root system $C_n$,
yields $W$-invariant eigenfunctions of the Sutherland differential
operator. The Sutherland differential operator after a suitable
gauge transformation becomes operator from the paper \cite{HO}. By
\cite{HO} $BC_n$ Jack polynomials are Weyl group invariant
eigenfunctions of this operator. Thus we obtain a representation
theoretic interpretation of Jack polynomials.

In \cite{HS} (see also \cite{St} for compact exposition of the
result and for its $q$-analog) a  one parameter family of the
$BC_n$ Jack polynomials was constructed by means of  spherical
functions on $G$. This family is a subfamily of the three
parameter  family of $BC_n$ Jack polynomials from
theorem~\ref{main} below. We also mention that in the
 one dimensional case
our construction reduces to the result of the paper \cite{Koel}.

 The structure of the paper is as follows. In Section 1,we give an
interpretation of the results on vector valued characters from
\cite{EK4} using the point of view of the theory of symmetric
spaces. Section 2 contains the main result: the construction of
Jack polynomials through  vector valued twisted spherical
functions on the symmetric space $G/K$. The proofs of the claims
from Section 2 are given in Section 3.  All constructions are
easier for the special case $m=n$,
$\varkappa_{(1)}+\varkappa_{(2)}=0$ and to gain a better
understanding,
 the reader
is advised to  consider this special case separately.

 The results of this paper can be  generalized to the case
of  quantum symmetric spaces \cite{ObSt}. The quantum version of
the construction yields the five parameter family of
Macdonald-Koornwinder polynomials, which are $q$ analog of Jack
polynomials.

{\bf Acknowledgments.} I would like to thank my adviser Pavel
Etingof for suggesting  the problem and  helpful explanations,
discussions and help with writing of  the text. I am  grateful to
Victor Ostrik and David Vogan for helpful discussions. I also
would like to thank Jasper Stokman who read carefully the
preliminary version of the text, corrected some  statements,
  gave me many useful references and helped
 me  make the
text much more readable. It was Jasper's suggestion to consider
the general case $m\ne n$  (the first version of the text
contained only the case $m=n$).

\section{Jack polynomials for the $A_{n-1}$ root system}
In this section we explain how to interpret the results of the
paper \cite{EK4} from the point of view of the theory of symmetric
spaces. This section does not contain any proofs. It just outlines
the results from \cite{EK4}.

\subsection{The space of spherical functions}
Let $K$ be the group $SL(n,\mathbb C)=SL(n)$ and $G=K\times K$.
The diagonal embedding $K\hookrightarrow K\times K$  gives rise to
the left and right action of $K$ on $G$. Let $V(\kappa)$ be the
representation  $S^{n\kappa}\mathbb C^n$ of $K$, $\kappa\in
\mathbb Z_+$.

The space of $K$-spherical functions $F_\kappa$ is defined by the
formula  $$ F_{\kappa}= \{f\in F(G, V(\kappa))| f(k g
k')=kf(g),\forall k,k'\in K, g\in G\}, $$ where $F(G,V(\kappa))$
is the space of $V(\kappa)$ valued polynomial functions on $G$.
The quotient space $G/K$ can be identified with the group $K$:
$(x,y)\to xy^{-1}$. Under this identification the right action of
group $K$ becomes the conjugation action and
 $$ F_\kappa=\{f\in
F(K,V(\kappa))|f(kk'k^{-1})=kf(k'),\forall k,k'\in K\}. $$

\subsection{The restriction of spherical functions to a maximal
torus and differential operators} Obviously the restriction of a
function from $F_{\kappa}$ to the
 maximal torus
$H=\{ e^{h(x)}|h(x)=diag(x_1,\dots,x_n), \mathop{\rm tr}
\nolimits (h(x))=0
 \}$ takes values in the one-di\-men\-si\-o\-nal space
$V(\kappa)[0]$, that is  it can be regarded  as a scalar function.

Furthermore, functions from $F_\kappa$
 are uniquely determined by their restriction to the maximal torus
(because the generic element of $K$ is conjugate to an element of
$H$). The Laplace operator on $K$ being restricted to $F_\kappa$
takes the  form: $$
\bar{L}_{\kappa}=\delta^{-1}\left(\sum_{i=1}^n\frac{\dd^2}{\dd
x_i^2}- \kappa(\kappa+1)\sum_{i<j}
\frac{1}{2\sinh^2((x_i-x_j)/2)}-\frac{n^3-n}{12}\right)\delta, $$
where $\delta=\prod_{i<j}\sinh((x_i-x_j)/2)$ and $x_i$
($i=1,\dots,n$) are the natural coordinates on $H$. The operator
$\delta\bar{L}_\kappa\delta^{-1}$ coincides, up to an additive
constant, with the Sutherland operator -- the
 Hamiltonian of the quantum $n$-body system
on the line with  interaction
 potential  $\kappa(\kappa+1)
\sinh^{-2}(y/2)$.

If $C_1,\dots, C_n$  are generators of the center $\mathfrak
Z(\mathfrak k)$ of the universal enveloping algebra $U(\mathfrak
k)$ of the  Lie algebra $\mathfrak k$ of the group $K$, then the
corresponding differential operators on $K$ map the  space
$F_\kappa$ into itself. Hence the restriction of such
 operator to the space $F_\kappa$ is a uniquely determined
 differential operator and it is possible to
 express it in terms of coordinates along  the maximal
 torus. Let us denote such differential operator on the torus by $R_{C_i}$.

\subsection{The connection with  quantum integrable systems}
Recall that an operator commuting with
the Hamiltonian of a quantum physical system is called a
{\it quantum integral} of this system.
Let $L$ be a differential operator in $n$ variables.
One says that $L$ defines a {\it completely integrable quantum Hamiltonian
system} if there exists a
set of $n$ algebraically independent quantum integrals
 $L_1,\dots, L_n$ which are differential operators and
 commute with each other. The collection
of operators  $L_1,\dots, L_n$
is called a complete system of quantum integrals for $L$.

It is clear from the above that we have following proposition

\begin{pro}
(1)\cite{OP} The Sutherland differential operator defines a
completely integrable system; moreover

(2)\cite{E} The operators $\delta R_{C_i}\delta^{-1}$,
$i=1,\dots,n$ form a complete system of quantum integrals for this
system.
\end{pro}

\subsection{Jack polynomials}
Below we give a definition of Jack polynomials following  the work
\cite{HO}. The action of the operator $\tilde{L}_\kappa=
\delta^{-\kappa}\bar{L}_\kappa\delta^\kappa$ on the space $\mathbb
{C}[P]^W$ of Weyl invariant Laurent polynomials on $H$ is
diagonalizable, and the eigenfunctions have the form
$J_\lambda^\kappa= m_\lambda+\sum_{\nu<\lambda}
s_{\lambda\nu}m_\nu$, where $m_\lambda(x)=\sum_{\nu\in W\lambda}
e^{(\nu,x)}$, where $\lambda\in P_+(SL(n))= \{\lambda\in\mathbb
Z^n|\lambda_1\ge \lambda_2\ge\dots \ge\lambda_n\}$ is an integral
dominant weight, defined up to shift,
 and $<$ is
the standard partial order  on weights. Polynomials
$J^\kappa_\lambda$ are called {\it Jack polynomials} for the root
system $A_{n-1}$.

\subsection{The formula for the spherical eigenfunction}
\label{normps} Let $L_\lambda$ be the finite dimensional
representation of $K$ with highest weight $\lambda\in P_+(SL(n))$
and $L^*_\lambda$ is its dual. Then the restriction of the finite
dimensional representation $L_\lambda\boxtimes L^*_\lambda$, of
$G=K\times K$  to the diagonally embedded subgroup $K$ yields the
representation $L_\lambda\otimes L^*_\lambda$ of $K$. The tensor
product $L_\lambda\otimes L^*_\lambda$ contains the vector
$w_0=\sum u_i\otimes u^*_i$ (where $\{u_i\}$ is a basis of
$L_\lambda$ and $\{ u^*_i\}$ is the corresponding dual basis of
$L^*_\lambda$), which is stable with respect to the (diagonal)
action of $K$. If $\lambda=\mu+ \kappa\rho$,
$\rho=\frac12\sum_{\alpha>0}\alpha$, $\mu\in P_+(SL(n))$,
 then $L_\lambda\otimes
L^*_\lambda$ also contains a unique copy of the
irreducible representation $V(\kappa)$.

 Thus for the tensor product one can write a decomposition
$$L_\lambda\otimes L^*_\lambda=V(\kappa)\oplus\oplus_{\mu\in
C_\lambda}L_\mu,$$ $C_\lambda\subset P_+(SL(n))$, $V(k)\ne L_\mu$
for all $\mu\in C_\lambda$ . Hence a vector $u\in L_\lambda\otimes
L^*_\lambda$ can be uniquely presented  in the form
$u=\tilde{u}+\sum_{\mu\in C_\lambda} u_\mu$, where $\tilde{u}\in
V(\kappa)$, $u_\mu\in L_\mu$. This decomposition gives an
embedding $s$: $V^*(\kappa)\hookrightarrow L_\lambda^* \otimes
L_\lambda$, $s(x)(u)=x(\tilde{u})$.

Let $v_1,\dots,v_N$ be a basis of representation $V(\kappa)$. Then
one can construct the function $$\Psi_\mu(g)=\sum v_i\langle
s(v_i^*),gw_0\rangle$$ (where $v^*_i$, $i=1,\dots,N$ is a basis
dual to the basis $v_i$, $i=1,\dots, N$). This function  belongs
to the space $F_\kappa$.  The maximal torus of $K$ is embedded in
$G$ by map $e^{h(x)}\mapsto (e^{h(x)},1)$, and the restriction of
$\Psi_\mu$ to this torus has the form
$$\Psi_\mu(x)=w_\kappa\langle s(w^*_\kappa), e^{h(x)}
w_0\rangle,$$ where $w^*_\kappa$, $w_\kappa$
 are the vectors such that $V(\kappa)^*[0]=
span\{ w^*_\kappa\}$, $V(\kappa)[0]=span\{ w_\kappa\}$ and
$\langle w_\kappa,w^*_\kappa\rangle=1$. Let us identify
$V(\kappa)[0]$ with $\mathbb C$ via $w_\kappa\to \langle w_\kappa,
v^*_\lambda\otimes v_\lambda\rangle$, where $v_\lambda$ and
$v_\lambda^*$ are the highest and lowest weight vectors for the
representations $L_\lambda$ and $L^*_\lambda$, respectively.

\subsection{The main theorem} An element of the center of the universal enveloping
algebra $U(\gotk\oplus \gotk)$
 acts on the space
$L_\lambda\boxtimes L^*_\lambda$ by  a constant.
In particular, the element $C_i\in U(\gotk)$,
embedded into
  $ U(\gotk\oplus \gotk)$ via $x\mapsto
(x,0)$ ($x\in\gotk$), does.  Hence $\Psi_\mu(x)$ is an
eigenfunction of $R_{C_i}$. Obviously, $\Psi_\mu(x)$ is
$W$-invariant when $\kappa$ is even and $W$ antiinvariant when $\kappa$ is odd.
 Moreover, the following theorem holds:

\begin{thrm}\cite{EK4} Under the above
identification of $V(\kappa)[0]$ with $\mathbb C$:

(1) $\Psi_0(x)=\delta(x)^\kappa$.

(2) $\Psi_\mu(x)$ is divisible by $\Psi_0(x)$ in algebra $\mathbb C[P]$ for
all $\mu\in P_+(SL(n))$.

(3) $\Psi_\mu(x)/\Psi_0(x)=J^\kappa_\mu(x)$.
\end{thrm}

\section{Jack polynomials for the $BC_n$ root system}
In this section we explain a representation theoretic construction
of the three parameter family of Jack polynomials for the root
system $BC_n$. We postpone all proofs (which are mostly lengthy
calculations) for the next section.
\subsection{The symmetric pair $(G,K)$ and restricted root
system}\label{restrr} In this section we introduce notations for
representation theoretic objects which are necessary for a further
exposition.

 Let $G$ be the
group $GL(m+n,\mathbb C)=GL(m+m)$
 and $\gotg$ its Lie algebra, where $m\ge n$. The conjugation by
 the diagonal matrix $J\in G$, $J_{ii}=1$, $i=1,\dots,m$, $J_{jj}=-1$,
 $j=m+1,\dots,m+n$ defines an involution $\Theta$:
 $\Theta(g)=JgJ^{-1}$. $\Theta$-invariant elements of $G$ form
 a subgroup $K=K_{(1)}\times
 K_{(2)}=GL(m)\times GL(n)\subset G$. The differential of this involution at
 the unit
 acts on the Lie algebra $\mathfrak g$ and induces the
  decomposition $\gotg=\gotk\oplus\mathfrak p$, $\gotk=\{x\in
\gotg|d\Theta(x)=x\}$, $\mathfrak p=\{x\in\gotg|d\Theta(x)=-x\}$ (
$\gotk$ is exactly the Lie algebra of the group $K$).The  subspace
$\mathfrak p\subset \gotg$ is not a subalgebra of $\gotg$ but
nevertheless  one can fix a maximal abelian subalgebra $\mathfrak
a$ inside $\mathfrak p$. We use notation $A$ for the corresponding
abelian subgroup. Below we fix  some particular subgroup $A$ which
we will work with.

Before fixing the choice of $A$ we introduce the Cartan subgroup
$H$ of the group $G$. The subgroup $H$ is conjugated to the
subgroup of the diagonal matrices by matrix $\tilde{J}$. Matrix
$\tilde{J}$ consists of the diagonal blocks:
$\tilde{J}_{ii}=\tilde{J}_{m+i,i}=
\tilde{J}_{i+m,i+m}=-\tilde{J}_{i,i+m}=\frac{1}{\sqrt{2}}$,
$i=1,\dots,n$, $\tilde{J}_{j+n,j+n}=1$, $j=1,\dots,m-n$ and all
other entries are zero:$$
\tilde{J}=\begin{pmatrix}\frac{I_n}{\sqrt{2}}&0&-\frac{I_n}{\sqrt{2}}\\
0&I_{m-n}&0\\ \frac{I_n}{\sqrt{2}}&0&\frac{I_n}{\sqrt{2}}
\end{pmatrix}.$$ We denote by the symbol $h(x,y,z)$  an element of
$H$ of the form $\tilde J e^{diag(x,y,z)}\tilde J^{-1}$, $x,z\in
\mathbb C^n, y \in\mathbb C^{m-n}$. Let $\mathfrak h$  be a Lie
algebra of $H$.

We put  $A$ to be equal to $\exp(\mathfrak p)\cap H$ or
$A=\{h(x,0,-x),x\in \mathbb C^n\}$. We use notation
$e^{a(x)}=h(x,0,-x)$, $x\in\mathbb C^n$ for the elements of $A$,
where $a(x)=\sum_{i=1}^n x_i(E_{i+m,i}-E_{i,i+m})\in\gotg$,
$E_{ij}$ is the notation for the $ij$-th matrix unit (basis in
$\gotg$).

The inclusion $a$: $\mathfrak a\hookrightarrow\mathfrak h$ induces
the projection
 $\mathfrak h^*\twoheadrightarrow \mathfrak a^*$. The root system
$R\subset \mathfrak h^*$ is mapped under this projection onto the
restricted root system. The restricted root system is isomorphic
to the root system $C_n$, in the case $n=m$, and to the root
system $BC_n$, in the case $m>n$. We use the notation $\Sigma$ for
the root system $BC_n$.   The short, medium and long positive
roots of $\Sigma$ are
 the
vectors:
\begin{align}
&\varepsilon_i& (1&\le i \le n),& \varepsilon_i&\pm\varepsilon_j&
(1&\le i<j\le n),& 2&\varepsilon_i& (1&\le i\le n),&
\end{align}
where $\varepsilon_i(a(x))=x_i$. The root multiplicities for the
short, medium and long  roots of $\Sigma$ are $t_1=2(m-n)$,
$t_2=2$, $t_3=1$. Below we use the half multiplicities
$s_i=t_i/2$, $i=1,2,3$.  The bilinear form $(\cdot,\cdot)$ on
$\mathfrak a^*$ is the standard one:
$(\varepsilon_i,\varepsilon_j)=\delta_{ij}$.

Let us also introduce  notations for a generalization of the Weyl
denominator:
\begin{equation*}\label{ta}
\delta_{p_1,p_2,p_3}(x)=
\prod_{\alpha\in\Sigma_+}\sinh^{p_\alpha}(\alpha(a(x))),
\end{equation*}
and for the vector: $$
\rho_{p_1,p_2,p_3}=\frac12\sum_{\alpha\in\Sigma_+}p_{\alpha}\alpha,
$$ where $p_{\pm\varepsilon_i}=p_1$, $p_{\pm\varepsilon_i\pm
\varepsilon_j}=p_2$ ($1\le i<j\le n$), $p_{\pm
2\varepsilon_i}=p_3$ ($1\le i<j\le n$). Below we always use the
convention that   for the given vector $\vec{p}\in \mathbb C^3$,
$p_\alpha$ means the same thing as in the previous formula if
$\alpha\in \Sigma$ and $p_\alpha=0$ if $\alpha\notin\Sigma$.

\subsection{Vector valued spherical functions} In this subsection
we define the space of $K$-equivariant twisted vector valued
functions on $G$. Functions from this space take values in the
space of the particular  representation of $K$. In principle one
has a lot of freedom in the choice of this representation and  it
is not clear why the chosen representation  is better than any
others. We will explain it in the next subsection.

First we define the space of "twisted" scalar valued functions on
$G/K$:  $$ \hat{F}= \{f\in F|f(gk)=f(g)
{\det}^{\varkappa_{(1)}}(k_{(1)}){\det}^{\varkappa_{(2)}}(k_{(2)}),
\forall k\in K, g\in G  \}, $$ where $F$ is the space of
polynomial functions on $G$ and $k=k_{(1)}k_{(2)}$, $k_{(i)}\in
K_{(i)}$. The desired space of spherical vector valued
$K$-spherical functions is a subspace of the tensor product
$\hat{F}\otimes U_{\vec{\varkappa}}$ where $U_{\vec{\varkappa}}$
is a representation of $K$ which we define below.

Let us fix notations for  finite dimensional representations of
$GL(r)$. A finite dimensional irreducible representation of
$GL(r)$ is encoded by its highest weight $\lambda\in\mathbb
Z^{r}$, $\lambda_1\ge \lambda_2\ge \dots\ge\lambda_r$, which is
dominant and integral. We denote the corresponding representation
by $L_\lambda$ and the set of highest weights by $P_+(GL(r))$. The
determinant representation in our notations is $L_{1^{r}}$,
$1^{r}$ is the $r$-dimensional vector consisting of ones,
$1^{r}=(1,\dots,1)$.

The representation $U_{\vec{\varkappa}}$ of the  group
$K=K_{(1)}\times K_{(2)}=GL(m)\times GL(n) $ is of the form
$U_{\vec{\varkappa}}=W_{\vec{\varkappa}}\boxtimes
V_{\vec{\varkappa}}$ and element $k=k_{(1)}k_{(2)}\in K$ acts on
$u=w\otimes v\in U_{\vec{\varkappa}}$ by formula
$k(u)=k_{(1)}(w)\otimes k_{(2)}(v)$.

The representation $W_{\vec{\varkappa}}$
 is the
an irreducible representation $L_{\mu}$,
$(\tilde{\varkappa}_{(1)}-\varkappa_{(1)})1^n+\varkappa_{(1)}1^m$
of the group $K_{(1)}=GL(m)$. Here
$\tilde{\varkappa}_{(1)},\varkappa_{(1)}$ are integers and $1^n\in
P_+(GL(m))$ is a vector of the form $(1,\dots,1,0,\dots,0)$ with
ones at first $n$ places and zeroes at other places. Why we choose
such representation  is clear from  lemma~\ref{1dW}.

The representation
 $V_{\vec{\varkappa}}$
 is an irreducible representation $L_{\lambda}\otimes\det^{\varkappa_{(2)}+
 \varkappa_{(1)}-\tilde{\varkappa}_{(1)}}$ of $K_{(2)}=GL(n)$,
 where
 $\lambda=(\varkappa_v(n-1),-\varkappa_v,\dots,-\varkappa_v).$
 More explicitly, the representation $L_\lambda$ is the representation
$S^{\varkappa_v n}\mathbb C^n$ (i.e. $\varkappa_vn$ symmetric
power of the vector representation) of $PGL(n)$ pulled back  to
$GL(n)$ i.e. the center acts trivially. The main reason why we
choose this representation is that zero weight space of
$L_\lambda$ is one-dimensional. Below we use notation
$\vec{\varkappa}=(\varkappa_{(1)},\varkappa_{(2)},\tilde{\varkappa}_{(1)},
\varkappa_v).$

\begin{rema}
In the simplest case, when $m=n$,  $W_{\vec{\varkappa}} \simeq
\det^{\tilde{\varkappa}_{(1)}}$  is one-dimensional. So
$U_{\vec{\varkappa}}$ becomes the representation
$\det^{\tilde{\varkappa}_{(1)}}\boxtimes
(\det^{-\tilde{\varkappa}_{(1)}}\otimes L_\lambda)$.
\end{rema}

Combining all components we get a definition of the space of
$K$-equivariant vector valued twisted spherical functions:
 $$
F_{\vec{\varkappa}}= \{f\in \tilde{F}\otimes
U_{\vec{\varkappa}}|f(kg)=kf(g), \forall k\in K, g\in G \}. $$ For
brevity we call the functions from this space spherical functions.

\begin{rema}
 In
the case $\varkappa_v=0$,
$\varkappa_{(1)}=\tilde{\varkappa}_{(1)}$, this space was studied
in the first part of the book \cite{HS}.
\end{rema}

\subsection{Properties of the spherical functions} In this
subsection we explain why the restriction of a spherical function
on the torus $A$ is  a scalar function.

Elementary arguments from the linear algebra show that the generic
 element $g$ of $G$ can be presented in the form $g=ke^{a(x)}k'$,
$k,k'\in K$, $x\in\mathbb C^n$ and this decomposition is unique up
to the action of the Weyl group.
  Because of the
 bi-$K$-equivariance of the functions from $F_{\vec{\varkappa}}$, any function
$f\in F_{\vec{\varkappa}}$ is uniquely determined by its
restriction to $A$.

The element $y$ of the group $M=Z_K(A)=GL(m-n)$ acts on a
spherical function $f$ following way:
$$yf(e^{a(x)})=f(ye^{a(x)})=f(e^{a(x)}y)=f(e^{a(x)})
\det(y)^{\varkappa_{(1)}}. $$

That is the restriction of a spherical function takes values in
the subspace $\tilde{W}_{\vec{\varkappa}}\boxtimes
V_{\vec{\varkappa}}$, where $\tilde{W}_{\vec{\varkappa}}$ is a
subspace on which $M$ acts by the character
$\det^{\varkappa_{(1)}}$.

\begin{lem}\label{1dW}
The subspace $\tilde{W}_{\vec{\varkappa}}$  is one-dimensional.
The subgroup $\tilde{K}_{(1)}=Z_{K_{(1)}}(M)=GL(n)$ acts on
$\tilde{W}_{\vec{\varkappa}}$ by the character
$\det^{\tilde{\varkappa}_{(1)}}$.
\end{lem}

The proof of the lemma is given at the section~\ref{prf}.

\begin{rema}
 In the simplest case  $m=n$,  the lemma is trivial,
 since then $W_{\vec{\varkappa}}$ reduces
to the one-dimensional representation
$\det^{\tilde{\varkappa}_{(1)}}$ of $K_{(1)}=\tilde{K}_{(1)}$.
\end{rema}

 Now let $T\subset K$  be the
subgroup $T=K\cap H=\{ k\in K|k=h(x,z,x), x\in\mathbb C^n,
z\in\mathbb C^{m-n}\}$. Observe that $T\subset
\tilde{K}_{(1)}\times M\times K_{(2)}$ and elements of $T$
commutes with $A$. An element $h_2=e^{diag(0,0,z)}$, $z\in\mathbb
C^n$ of the Cartan subgroup of $K_{(2)}$ can be presented in the
form $h_2=t h_1$, where $h_1=e^{diag(-z,0,0)}\in \tilde{K}_{(1)}$
and $t\in T$. Hence using $\det(h_1)\det(h_2)=1$ one gets for a
spherical function $f$: $$
h_2f(e^{a(x)})=f(h_1te^{a(x)})=h_1f(e^{a(x)})t=
\det(h_2)^{\varkappa_{(1)}+\varkappa_{(2)}-\tilde{\varkappa}_{(1)}}f(e^{a(x)}).$$

 That means that the restriction of $f$ to $A$ takes values in
the one-dimensional space $\tilde{W}_{\vec{\varkappa}}\boxtimes
\tilde{V}_{\vec{\varkappa}}$, with $\tilde{V}_{\vec{\varkappa}}=
V_{\vec{\varkappa}}[(\varkappa_{(1)}+\varkappa_{(2)}-\tilde{\varkappa}_{(1)})1^n]
\simeq \mathbb C$.
 That is,
 the restriction of a spherical function $f$ to $A$ is a scalar
valued function.

\subsection{The center of $U(\mathfrak g)$ and radial parts of biinvariant
differential operators} In this subsection we explain the
correspondence between elements of the center of the universal
enveloping algebra and the Weyl group invariant differential
operators on $A$. This correspondence is given by the radial parts
of the differential operators. Finally we calculate the radial
part of the  Casimir operator.

The universal enveloping algebra $U(\gotg)$ may be identified with
the algebra of the left $G$-invariant differential operators on
$G$. Namely, the element $x\in\gotg$ gives the differential
operator $D_{x}f(g)=\frac{d}{dt}f(ge^{tx})|_{t=0}$ and this map
can be extended to $U(\gotg)$: $D_{xy}f=D_{x}D_yf$, for $x,y\in
U(\gotg)$. A differential operator corresponding to an element of
the center $\mathfrak Z(\gotg )$  of the universal enveloping
algebra $U(\gotg)$ is bi-$G$-invariant, hence it preserves the
 space
$F_{\vec{ \varkappa}}$. As any function from $F_{\vec{\varkappa}}$
is uniquely determined by its  restriction to $A$, the
differential operator $D_C$, $C\in \mathfrak{ Z(g)}$ can be
written in terms of coordinates along $A$, and the resulting
operator is a differential operator with coefficients in
$End(\tilde{U}_{\vec{\varkappa}})=\mathbb C$. We call this
expression the {\it radial part of} $D_C$ on $F_{\vec{\varkappa}}$
 and denote it $R_C$.

\begin{rema}In reality, not any function
on the torus $A$ is the restriction of  a
$K$-bi\-equi\-va\-ri\-ant function. Later we will show that the
restriction of the  functions from $F_{\vec{\varkappa}}$ span the
space of Laurent polynomials of $e^{2x_i}$ ($x_i$, $i=1,\dots,n$,
are the coordinates along the torus $A$) satisfying vanishing
conditions at zero locus locus of Weyl determinant (see
lemma~\ref{dd}), on which $W$ acts by the character $\chi$ (see
lemma~\ref{Wact}). Here $W$ is the Weyl group for the $BC_n$ root
system and $\chi$ is the $\mathbb Z_2$-character of $W$. But the
standard reasoning (see for example \cite{HS} page 16) shows that
a differential operator is uniquely determined by its action on
the space of $\chi$-$W$ invariant Laurent
 polynomials, hence the radial part of $D_C$ is uniquely defined.
\end{rema}

The center $\mathfrak Z(\gotg)\subset U(\gotg)$ contains the Casimir
 $C_2=
\sum_{1\le i,j\le n+m} E_{ij}E_{ji}$. The radial part $R_{C_2}$
can be calculated explicitly. Below we use following notations
$\vec{\kappa}=(\kappa_1,\kappa_2,\kappa_3)=(\varkappa_{(2)}-\varkappa_{(1)},
\varkappa_v,\tilde{\varkappa}_{(1)}-\varkappa_{(2)})$,
$\vec{r}=\vec{\kappa}+\vec{s}$. We consider only the case when
$\vec{\kappa}$ satisfies the following condition:
\begin{equation}\label{restrk}
\kappa_3\ge \kappa_1+\kappa_3\ge 0.
\end{equation}

\begin{thm}\label{CM}
The  second order differential operator $R_{C_2}$ has the form
\begin{multline}\label{CMf}
2(R_{C_2}\psi)(x)=\delta^{-1}_{\vec{s}}(x)\left(
\Delta_{A}-u_{\vec{r}}(x)
+C_{\vec{\varkappa}}\right)\delta_{\vec{s}}(x)\psi(x),\\
C_{\vec{\varkappa}}=\frac{m+n-(m+n)^3}{6}+\frac{(m-n)^3-m+n}{6}+(\varkappa_{(1)}+
\varkappa_{(2)})^2n+2(m-n)\varkappa_{(1)}^2,
\end{multline}
\begin{equation}
u_{\vec{r}}(x)= \sum_{\alpha\in\Sigma_+}\frac{{r}_\alpha
({r}_{\alpha}+2{r}_{2\alpha}-1)(\alpha,\alpha)}{{\sinh(\alpha(a(x)))}^2},
\end{equation}
where   $\Delta_{A}=\sum_{i=1}^n\frac{\dd^2}{\dd x_i^2}$
 is the Laplace operator.
\end{thm}

 The operator $L^{CM}_{\vec{r}}=
\Delta_{A}-u_{\vec{r}}(x)$ is called the Calogero-Moser operator
for the root system $BC_n$.

\begin{rema} When $\varkappa_{(1)}+\varkappa_{(2)}=0$,  $m=n$,
 the coefficient
$C_{\vec{\varkappa}}$ is equal to $-2(\rho,\rho)$ ($\rho$ is half
sum of the positive roots of the root system $A_{m+n-1}$).
\end{rema}

\subsection{The spherical representations of $G$}
In this subsection we describe the finite dimensional
representations of $G$ containing the representation
$U_{\vec{\varkappa}}$ of $K$. Proofs of the statements use the
Littlewood-Richardson rule and are given at the next section.

First we introduce some new definitions and notations. Let us
denote by $P^{\varkappa_{(1)},\varkappa_{(2)}}$,
$\varkappa_{(1)},\varkappa_{(2)}\in \mathbb Z$,
($\varkappa_{(1)}\ge\varkappa_{(2)}$ by (\ref{restrk}))  the
subset of $P_+(GL(n+m))$ consisting of $\lambda\in P_+(GL(n+m))$
such that $G$-representation $L_\lambda|_K$ contains a copy of
$\det^{\varkappa_{(1)}}\otimes\det^{\varkappa_{(2)}}$.

\begin{lem}\label{S}
 $\lambda\in P^{\varkappa_{(1)},\varkappa_{(2)}}$ if and only if
\begin{align*}
\lambda_{j}+\lambda_{m+n+1-j}&=\varkappa_{(1)}+\varkappa_{(2)}&
(j&=1,\dots,n)&\\ \lambda_{n+j}&=\varkappa_{(1)}&
(j&=1,\dots,m-n)&\\ \lambda_n&\ge \varkappa_{(1)}.&
\end{align*}
Moreover, if $\lambda\in P^{\varkappa_{(1)},\varkappa_{(2)}}$ then
$L_\lambda|_K$ contains a unique copy of the representation
${\det}^{\varkappa_{(1)}}\boxtimes{\det}^{\varkappa_{(2)}}$.
\end{lem}

\begin{rema}
In the case when all $\varkappa_{(i)}$ are zero the last lemma
follows from the fact that $(G,K)$ is a symmetric pair (see
\cite{He2} Chapter V, Theorem 4.1).
\end{rema}

We denote by the symbol $P_+^{BC}$ the set of $n$-tuples of
non-negative integers $\mu\in \mathbb Z^n_+$ such that
$\mu_1\ge\mu_2\ge\dots\ge\mu_n\ge 0$.  Now consider the map
$\tau$: $P_+^{BC}\to \mathbb Z^{m+n}$, where $$ \tau(\mu)=
(\mu_1+\hat{\varkappa},\dots,
\mu_n+\hat{\varkappa},\varkappa_{(1)},\dots,\varkappa_{(1)},
\hat{\varkappa}-\mu_n,\dots,\hat{\varkappa}-\mu_1),$$
$\frac{\varkappa_{(1)}+\varkappa_{(2)}}2=\hat{\varkappa}.$
   Lemma~\ref{S} says that $\lambda\in
P^{\varkappa_{(1)},\varkappa_{(2)}}$ if and only if
$\lambda=\tau(\rho_{-\kappa_1,0,0}+\mu)$ for some $\mu\in
P_+^{BC}$. That is, lemma~\ref{S} implies that the set
$P^{\varkappa_{(1)},\varkappa_{(2)}}$ is
 isomorphic to $P_+^{BC}$.

\begin{defi} A finite dimensional irreducible representation
 $L_\lambda$ of $G$ is called {\it $\vec{\varkappa}$-spherical}
($\vec{\varkappa}$ is integer and satisfies (\ref{restrk})), if
$\lambda\in P^{\varkappa_{(1)},\varkappa_{(2)}}$ and
$L_\lambda|_K\supset U_{\vec{\varkappa}}$. Let us denote by
$P^{\vec{\varkappa}}$ a set of weights $\lambda$ such that
$L_\lambda$ is $\vec{\varkappa}$-spherical.
\end{defi}

\begin{lem}\label{s1} If $\mu\in P_+^{BC}$ then
$\lambda=\tau(\mu+\rho_{\vec{\kappa}})\in P^{\vec{\varkappa}}$,
and in this case $L_\lambda|_K$ contains a unique copy of the
representation $U_{\vec{\varkappa}}$.
\end{lem}

Also lemma~\ref{S} and (\ref{restrk}) imply that
$\tau(\mu+\rho_{\vec{\kappa}})\in
P^{\varkappa_{(1)},\varkappa_{(2)}}$ for $\mu\in P^{BC}_+$.

\begin{rema}The  statement converse to the lemma~\ref{s1} also holds.
It follows
from the statements that are stated below (see
corollary~\ref{degr}). Let us also remark that in the simplest
case ( $n=m$, $\varkappa_{(1)}+\varkappa_{(2)}=0$) the map $\tau$
does not depend
 on $\vec{\varkappa}$.
\end{rema}

\subsection{Spherical functions through $\vec{\varkappa}$-spherical
 representations}\label{eig}
 For $\lambda\in P^{\vec{\varkappa}}$ we have a
decomposition $$L_\lambda|_{K}=
U_{\vec{\varkappa}}\oplus_{(\mu,\mu')\in C_\lambda}L_\mu\boxtimes
L_{\mu'},$$ $C_\lambda\subset P_+(GL(m))\oplus P_+(GL(n))$. Hence
a vector $v\in L_\lambda$ can be uniquely presented in the form
$v=\tilde{v}+\sum_{(\mu,\mu')\in C_\lambda} v_{\mu,\mu'}$,
$\tilde{v}\in U_{\vec{\varkappa}}$, $v_{\mu,\mu'}\in
L_\mu\boxtimes L_{\mu'}$. It allows us  to define an embedding
$s$: $U^*_{\vec{\varkappa}}\hookrightarrow L^*_\lambda|_K$,
$s(x)(v)=x(\tilde{v})$.

Let $u_{\varkappa_{(1)},\varkappa_{(2)}}(\lambda)\in L_\lambda$,
$\lambda\in P^{\varkappa_{(1)},\varkappa_{(2)}}$
 be a vector  such that $$span\{ u_{\varkappa_{(1)},\varkappa_{(2)}}(\lambda)
\}={\det}^{\varkappa_{(1)}}\boxtimes{\det}^{\varkappa_{(2)}}.$$
The vector $u_{\varkappa_{(1)},\varkappa_{(2)}}(\lambda)$ is
unique up to normalization. Now consider the function
$\Psi_{\mu}$: $G\to U_{\vec{\varkappa}}$, $\mu\in P_+^{BC}$, given
by $$\Psi_\mu(g)= \sum v_i\langle s(v_i^*),
gu_{\varkappa_{(1)},\varkappa_{(2)}}(\lambda)\rangle,$$ where
$\lambda=\tau (\mu+\rho_{\vec{\varkappa}})$
 and
$v_i$, $i=1,\dots,N$ is a basis of the representation
$U_{\vec{\varkappa}}\subset L_\lambda|_K$ and $v^*_i$,
$i=1,\dots,N$ is a dual basis. It is easy to see that $\Psi_\mu\in
F_{\vec{\varkappa}}$ and that its restriction to $A$ is equal to
$w_{\vec{\varkappa}}\langle s( w^*_{\vec{\varkappa}}),
e^{a(x)}u_{\varkappa_{(1)},\varkappa_{(2)}}(\lambda)\rangle$,
where $span\{ w_{\vec{\varkappa}}\}= \tilde{U}_{\vec{\varkappa}}$,
$span\{ w^*_{\vec{\varkappa}}\}= \tilde{U}^*_{\vec{\varkappa}}$
and $\langle w^*_{\vec{\varkappa}}, w_{\vec{\varkappa}}\rangle=1$.
We identify $\tilde{U}_{\vec{\varkappa}}$ with $\mathbb C$ via
 $w_{\vec{\varkappa}}\to 1$ and to simplify
notations we write $$\Psi_\mu(x)=\langle w^*_{\vec{\varkappa}},
e^{a(x)}u_{\varkappa_{(1)},\varkappa_{(2)}}(\lambda)\rangle.$$

\begin{rema}Such definition for $\Psi_\mu$ has a flaw. The vectors
$w^*_{\vec{\varkappa}}$, $u_{\varkappa_{(1)},\varkappa_{(2)}}$ are
determined up to multiplication on a
 constant. Hence $\Psi_\mu$ is also defined up to multiplication by a constant.
We fix this constant at the end of  section~\ref{last}.
\end{rema}

\subsection{Eigenvalues of the radial parts of biinvariant
  differential operators}
The constructed function is an eigenfunction of some collection
 of operators.
Indeed, the  elements $C_r= \sum_{1\le i_1,\dots,i_r\le m+n}
E_{i_1i_2} E_{i_2i_3}\dots E_{i_ri_1}$ generate the center
$\mathfrak Z(\gotg)$ of $U(\gotg)$. By the Harish-Chandra theorem
we have \begin{gather*}
C_r|_{L_\lambda}=\sum_{j=1}^{n+m}(\lambda_j^r+\mbox{terms of lower
degree on }\lambda) Id_{L_\lambda},\\ C_2|_{L_\lambda}=
(\lambda,\lambda+2\rho)Id_{L_\lambda}, \end{gather*} where
$\rho=\frac12(m+n-1,m+n-3,\dots,-m-n+1)$.

\begin{rema} Actually, by induction on $r$ one can prove more
precise formula:
$$C_r|_{L_\lambda}=\sum_{j=1}^{n+m}((\lambda_j+\rho_j)^r-\rho_j^r)
Id_{L_\lambda},$$ but for our purposes the weaker formula is
sufficient. \end{rema}

 Using formulas from lemma~\ref{S} and the
previous formula, for any positive $j$ one gets:
\begin{equation}\label{c2}
R_{C_2}\Psi_\mu=(2(\mu+\rho_{\vec{r}},
\mu+\rho_{\vec{r}})+\frac{C_{\vec{\varkappa}}}2)\Psi_\mu
\end{equation}
\begin{equation}\label{Cch}
R_{C_{2j}}\Psi_\mu=2(\sum_{i=1}^n\mu_i^{2j}+
\mbox{terms of lower degree})\Psi_\mu,
\end{equation}
\begin{equation}\label{nech}
R_{C_{2j+1}}\Psi_\mu=
((2r+1)(\varkappa_{(1)}+\varkappa_{(2)})\sum_{i=1}^n\mu_i^{2j}+\mbox{terms
of lower degree})\Psi_\mu.
\end{equation}

\begin{rema}
 At the simplest case, when $n=m$, $\varkappa_{(1)}+ \varkappa_{(2)}=0$, formulas
(\ref{c2})-(\ref{nech}) are  simpler: $$ R_{C_{2j}}\Psi_\mu
=2(\sum_{i=1}^n(\mu+\rho_{\vec{r}})_i^{2j}-(\rho_{\vec{s}})_i^{2j})\Psi_\mu,
$$ $$ R_{C_{2j+1}}\Psi_{\mu}=0, $$
  for any positive integer $j$.
\end{rema}

\subsection{The Weyl group invariance and factorization of the
spherical function} The Weyl group $W$ of the $BC_n$ root system
 naturally maps  onto the group
$S_n$. We denote this map by $q$. The  group $W$ has two
independent $\mathbb Z_2$-characters. Indeed if $t: W\to
GL(n,\mathbb Z)$ is the tautological representation of this group,
then $\chi_0(w)=\det(t(w))$ and $\chi_1(w)=(-1)^{q(w)}$ for a
basis of $\mathbb Z_2$-characters.  For the character
$\chi=\chi_0^{\kappa_1+\kappa_3}\chi_1^{\kappa_1+\kappa_2+\kappa_3}$
the following statement holds:

\begin{lem}\label{Wact}
 A function $f(e^{a(x)})\in F_{\vec{\varkappa}}$ transforms under
the action of $W$ by the character $\chi$. Besides, we have
$f(e^{a(x+\pi i e_j)})= (-1)^{\kappa_1}f(e^{a(x)})$ for
$j=1,\dots,n$.
\end{lem}

The function $\Psi_\mu(x)$ is a Laurent polynomial of $e^{x_l}$
because $L_\lambda$ is a  polynomial representation. Let us denote
the space of Laurent polynomials in $e^{x_l}$, $l=1,\dots, n$
 by $\mathbb C[P]$.

\begin{lem}\label{dd} Any $f(e^{a(x)})$, $f\in F_{\vec{\varkappa}}$
 is divisible by
$\delta_{\vec{\kappa}}$ in the algebra $\mathbb C[P]$.
\end{lem}

Lemma~\ref{Wact} and lemma~\ref{dd}
 imply that the function $\Psi_\mu(x)/\delta_{\vec{\kappa}}(x)$
belongs to the space $\mathbb C[P]^W$ of $W$-invariant Laurent
polynomials in $e^{2x_l}$, $l=1,\dots,n$.
 Moreover, the following corollary holds

\begin{cor}\label{degr}
\begin{enumerate}
\item $\Psi_\mu(x)/\delta_{\vec{\kappa}}(x)
=\sum_{\nu\le\mu} d_{\mu\nu}
m_\nu(x)$, where
$m_\nu$ is the orbitsum $m_\nu=\sum_{\lambda\in
W\nu}e^{2(\lambda,x)}$ and $\le$ is
the standard $BC_n$ dominance order.
\item If $\lambda\in P^{\vec{\varkappa}}$ then
$\lambda=\tau(\mu+\rho_{\vec{\varkappa}})$, where $\mu\in
P_+^{BC}$.
\end{enumerate}
\end{cor}

\subsection{The definition and properties of  Jack polynomials}

 Consider the  operator
$\tilde{L}_{\vec{r}}=\delta^{-1}_{\vec{r}} L^{CM}_{\vec{r}}
\delta_{\vec{r}}$. For this operator the following proposition
holds.

\begin{pro} \cite{HO}\label{HOt} The operator $\tilde{L}_{\vec{r}}$
maps the  space $\mathbb C[P]^W$ of Weyl group invariant Laurent
polynomials into itself. Moreover, it is triangular with respect
to  the basis of the orbitsums $m_\lambda(x)=\sum_{\nu\in
W\lambda} e^{2(\nu,x)}$ i.e.
 $$\tilde{L}_{\vec{r}}
 m_\mu=4(\mu+\rho_{\vec{r}},\mu+\rho_{\vec{r}})m_\mu
+\sum_{\nu<\mu}\alpha_{\mu\nu}m_\nu.$$
\end{pro}

This proposition implies that one can uniquely determine the
Laurent polynomial $J^{\vec{r}}_\mu=m_\mu+
\sum_{\nu<\mu}s_{\mu,\nu}m_\nu$ by the condition
$\tilde{L}_{\vec{r}}J^{\vec{r}}_\mu=
4(\mu+\rho_{\vec{r}},\mu+\rho_{\vec{r}})J^{\vec{r}}_\mu$. Indeed,
$\rho_{\vec{r}}$ is a dominant weight, hence $\nu<\mu$ implies
$(\nu+\rho,\nu+\rho)<(\mu+\rho,\mu+\rho)$. Thus the operator
$\tilde{L}_{\vec{r}}$ being restricted to the finite dimensional
space $span\{m_\nu,\nu\ge\mu\}$ is diagonalizable with the
distinct eigenvalues. Hence $J^{\vec{r}}_\mu$ is uniquely
determined. The polynomials $J^{\vec{r}}_\mu$
 are called {\it Jack polynomials} for the $BC_n$ root system.

It is easy to see that the operator $L^{CM}_{\vec{r}}$ is
self-adjoint with respect to the standard inner product
$(f,g)=\int_{A^*} f(x) \overline{g(x)}dx$ (here the bar means
complex conjugation and $A^*=\{ e^{a(x)}|{\rm Re} x=0\}$). This
fact implies that $J_\mu^{\vec{r}}$ is orthogonal to
$J_\nu^{\vec{r}}$ if $(\mu+\rho_{\vec{r}},\mu+\rho_{\vec{r}})\ne
(\nu+\rho_{\vec{r}},\nu+\rho_{\vec{r}})$. In fact,
 an even stronger statement holds:
\begin{pro}\cite{HO}
\label{ort} The Jack polynomials $J^{\vec{r}}_\mu$, $\mu\in
P_+(Sp(n))$ form an orthogonal basis in the space $\mathbb
C[P]^W$. That is, $$(J^{\vec{r}}_\mu,J^{\vec{r}}_\nu)_{\vec{r}}=
\int_{A^*} \delta_{\vec{r}}(x)\overline{\delta_{\vec{r}}(x)}
J^{\vec{r}}_\mu \overline{J^{\vec{r}}_\nu} dx=0,$$ if $\mu\ne
\nu$.
\end{pro}

This proposition and theorem~\ref{CM} imply
\begin{cor}\label{dne0}
 The coefficient $d_{\mu\mu}$ at the expansion
$\Psi_\mu/\delta_{\vec{\kappa}}=\sum_{\nu\le \mu}d_{\mu\nu}m_\nu$ is
not zero.
\end{cor}

\subsection{The formulation of the main result}\label{last} The last
corollary and lemma~\ref{as} (see next section) imply $\langle
v_\lambda^*,u_{\varkappa_{(1)},\varkappa_{(2)}}(\lambda)\rangle\ne
0$, $\langle v_\lambda,w_{\vec{\varkappa}}\rangle\ne 0$. Let us
renormalize the function $\Psi$: $$ \tilde{\Psi}_\mu=\frac{\langle
s(w^*_{\vec{\varkappa}}), e^{a(x)}
u_{\varkappa_{(1)},\varkappa_{(2)}}(\lambda) \rangle}{\langle
s(w^*_{\vec{\varkappa}}),v_\lambda\rangle \langle v^*_\lambda,
u_{\varkappa_{(1)},\varkappa_{(2)}}(\lambda) \rangle}, $$ where
$\lambda=\tau(\mu+\rho_{\vec{\varkappa}})$ and $v_\lambda$,
$v^*_\lambda$, $\langle v_\lambda,v^*_\lambda\rangle=1$ are the
highest and lowest weight vectors for the $G$-representations
$L_\lambda$ and $L^*_\lambda$, respectively. Now, the  function
$\tilde{\Psi}_\mu$ does not depend on the choice either
$w^*_{\vec{\varkappa}}$ or
$u_{\varkappa_{(1)},\varkappa_{(2)}}(\lambda)$ (see the discussion
at the end of section~\ref{eig}).

The following theorem explains how to get  Jack polynomials from
the
 spherical functions. It also gives some details about the radial parts
of $C_r$, $r\in\mathbb N$.

\begin{thm}\label{main}
\begin{enumerate}
\item $\tilde{\Psi}_\mu/\delta_{\vec{\kappa}}=J^{\vec{r}}_\mu$.
\item The radial parts $R_{C_{2i}}$, $i\in\mathbb N$
are pairwise commutative differential operators in $n$ variables
of the form $$ R_{C_{2i}}=2^{1-2i}\sum_{j=1}^n \frac{\dd^{2i}}{\dd
x_j^{2i}}+ \sum_{J,|J|<2i}a_j(x)\frac{\dd^J}{\dd x^J}.$$
\end{enumerate}
\end{thm}
\begin{rema} In the case $m=n$, $\varkappa_{(1)}+\varkappa_{(2)}=0$,
the radial parts $R_{C_{2i+1}}$, $i\in \mathbb N$ are zero. In the
general case, the radial parts $R_{C_{2i+1}}$, $i\in \mathbb N$
can be expressed through $R_{C_{2j}}$.
\end{rema}

The second item of the theorem implies the complete integrability
(see previous section for the definition) of the quantum
Hamiltonian system defined by the Calogero-Moser operator
$L^{CM}_{\vec{r}}$. The first proof of the complete integrability
of this system was given by Olshanetsky and Perelomov \cite{OP}.
The quantum integrals $R_{C_{2i}}$ from the second part the
theorem coincide with the integrals from the paper
\cite{OP},because after conjugation by $\delta_{\vec{\kappa}}$
they are diagonal in the basis of Jack polynomials, with the same
eigenvalues as operators from \cite{OP}.

\section{Proofs}\label{prf}
We consider  an element $g\in G$ as a $3\times 3$  block matrix:
$$ g=\begin{pmatrix} g^{11}&g^{12}&g^{13}\\ g^{21}&g^{22}&g^{23}\\
g^{31}&g^{32}&g^{33}
\end{pmatrix},
$$ in which the $11$-th, $13$-th, $31$-th and $33$-th blocks are
$n\times n$ matrices, $12$-th and $21$-th block are $n\times
(m-n)$ and $ (m-n)\times n$ matrices, respectively and $22$-th
block is a $(m-n)\times(m-n)$ matrix. We denote
 these blocks by $g^{ij}$, $i,j=1,2,3$, and the matrix elements of $g$ by
$g^{ij}_{st}$, $i,j=1,2,3$, $s=1,\dots,n-(1+(-1)^i)\frac{m}2$,
$t=1,\dots,n-(1+(-1)^j)\frac{m}2$.

\subsection{Calculation of the radial part for
the Casimir element}\begin{proof}[Proof of theorem\ref{CM}] Using
the formulas $E^{kl}_{ij}E^{k'l'}_{ij}=0$,
$e^{sE_{ij}^{kl}}=1+sE_{ij}^{kl}$ for $i\ne j$ one gets :
\begin{multline*}
e^{t(\sinh x_j\cosh x_i E^{11}_{ij}+\sinh x_i\cosh
x_jE^{33}_{ij})} e^{a(x)} e^{sE^{31}_{ji}}
e^{st\sinh(x_i+x_j)\sinh(x_i-x_j)(E_{jj}^{33}- E^{11}_{ii})}\\
\times e^{st(\sinh x_j\cosh x_j E^{31}_{jj}- \sinh x_i\cosh x_i
E^{31}_{ii})} =e^{a(x)}e^{sE^{31}_{ji}} e^{t
\sinh(x_j-x_i)\sinh(x_i+x_j)E^{13}_{ij}}\\ \times e^{t(\sinh
x_j\cosh x_j E^{11}_{ij}+ \sinh x_i\cosh x_i E^{33}_{ij})}+O(t^2)+
O(s^2).
\end{multline*}
Substituting RHS and LHS of the last equation into the argument of
 a function $f\in F_{\vec{\varkappa}}$  and
 taking the derivative $\frac{d}{dt}$ at the point $t=0$, one gets
\begin{multline}\label{ff}
(\sinh x_j\cosh x_i E^{11}_{ij}+\sinh x_i\cosh x_j E^{33}_{ij})
f(e^{a(x)}e^{s E^{31}_{ji}})+sf(e^{a(x)}e^{sE^{31}_{ji}})\\
\times\sinh(x_j-x_i) \sinh(x_i+x_j)(E^{33}_{jj}-E^{11}_{ii})+
s(\sinh x_j\cosh x_j D_{E^{31}_{jj}}f(e^{a(x)}e^{sE^{31}_{ji}})\\-
\sinh x_i\cosh x_i D_{E^{31}_{ii}} f(e^{a(x)}e^{sE^{31}_{ji}}))=
\sinh(x_i-x_j)\sinh(x_i+x_j)\\
\times D_{E^{13}_{ij}}f(e^{a(x)}e^{sE^{31}_{ji}}).
\end{multline}
Substituting $s=0$ to (\ref{ff}) and changing $i$ and $j$, we have
\begin{multline}\label{eij}
(\cosh x_j\sinh x_i E^{33}_{ji}+\sinh x_j\cosh x_i E^{11}_{ji}) f(e^{a(x)})\\=
\sinh(x_i-x_j)\sinh(x_i+x_j)D_{E^{31}_{ji}}f(e^{a(x)}).
\end{multline}
 Taking the derivative $\frac{d}{ds}$
of formula (\ref{ff}) at the point $s=0$  and using (\ref{eij})
yields
\begin{multline}\label{Dij}
D_{E_{ji}^{31}E_{ij}^{13}} f(e^{a(x)})=f(e^{a(x)})(E^{33}_{jj}-E^{11}_{ii})
-(\sinh x_j\cosh x_i E^{11}_{ij}+\sinh x_i\\ \times\cosh x_j
E^{33}_{ij})
\frac{(\cosh x_j\sinh x_i E^{33}_{ij} +\sinh x_j \cosh x_i E^{11}_{ji})}
{\sinh^2(x_j-x_i)\sinh^2(x_i+x_j)}f(e^{a(x)})\\
+\frac{\sinh x_i\cosh x_i D_{E^{31}_{ii}}-\sinh x_j\cosh x_j D_{E^{31}_{jj}}}
{\sinh(x_i-x_j)\sinh(x_i-x_j)}f(e^{a(x)}).
\end{multline}
The elements $e^{a(x)}\in G$, $x\in\mathbb C^n$ form a commutative
subgroup isomorphic to  an $n$-dimensional torus
($e^{a(x+y)}=e^{a(x)}e^{a(y)}$), hence
\begin{equation}\label{Dii}
\frac{\dd f}{\dd
x_i}(e^{a(x)})=(D_{E_{ii}^{31}}-D_{E_{ii}^{13}})f(e^{a(x)}).
\end{equation}
 Substituting the formulas for the right action of $\mathfrak k$ on
the space $\tilde{W}_{\vec{\varkappa}}\boxtimes
\tilde{V}_{\vec{\varkappa}}$:
$E^{33}_{ij}E^{33}_{ji}=\varkappa_v(\varkappa_v+1)$,
$E^{11}_{ij}=0$, $E^{33}_{ii}=\varkappa_{(1)}+\varkappa_{(2)}-
\tilde{\varkappa}_{(1)}$ $E^{11}_{ii}=\tilde{\varkappa}_{(1)}$,
$i\ne j$ and for the left action $E^{11}_{ii}=\varkappa_{(1)}$,
$E^{33}_{ii}=\varkappa_{(2)}$, $E^{11}_{ij}=E^{33}_{ij}=0$, $i\ne
j$ to the formula (\ref{Dij}) and using (\ref{Dii}) one gets:
\begin{multline}\label{DDij}
D_{E^{31}_{ij}E^{13}_{ji}}+D_{E^{31}_{ji}E^{13}_{ij}}+
D_{E^{13}_{ij}E^{31}_{ji}}+D_{E^{13}_{ji}E^{31}_{ij}}\\=
-\varkappa_v(\varkappa_v+1)\left(\frac{1}{\sinh^2(x_i+x_j)}
+\frac{1}{\sinh(x_i-x_j)}\right)\\+
\frac{1}{\sinh(x_i-x_j)\sinh(x_i+x_j)} \left(\sinh
2x_i\frac{\dd}{\dd x_i} -\sinh 2x_j\frac{\dd}{\dd x_j}\right).
\end{multline}

The calculation of $D_{E^{31}_{ii}E^{13}_{ii}}$  for $m\ge n>1$ is
absolutely the same as in the case $m=n=1$. We make this
calculation for $n=1$.

 For $f\in F_{\vec{\varkappa}}$ the following equation holds
$$ f(z)=
\left(\frac{z^{13}}{z^{31}}\right)^{\frac{\tilde{\varkappa}_{(1)}-
\varkappa_{(1)}}2}
\left(\frac{z^{11}}{z^{33}}\right)^{\varkappa_{(1)}+\frac{\varkappa_{(2)}
+\tilde{\varkappa}_{(1)}}2}
{\det(z)}^{\frac{\varkappa_{(1)}+\varkappa_{(2)}}2} f(e^{a(x)}),
$$ where $x=\mathop{\rm arcsinh}\nolimits
\left(\sqrt{\frac{z^{13}z^{31}}{\det(z)}}\right)$. Hence we have
\begin{multline}\label{DDii}
f(e^{a(x)}e^{sE^{31}}e^{tE^{13}})= \left(\frac{\sinh x+t\cosh x+st
\sinh x} {\sinh x+s\cosh
x}\right)^{\frac{\tilde{\varkappa}_{(1)}-\varkappa_{(1)}}2}\\\times
\left(\frac{\cosh x+s\sinh x} {\cosh x+t\sinh x+st\cosh
x}\right)^{\varkappa_{(1)}+\frac{\varkappa_{(2)}+\tilde{\varkappa}_{(1)}}2}f(e^{a(y)}),
\end{multline}
where $y=\mathop{\rm arcsinh}\nolimits(\sinh( x)\sqrt{(1+s \coth
x)((1+st)+t \coth x)}))$. Taking the derivative $\frac{\dd^2}{\dd
s\dd t}$ of (\ref{DDii}) at the point $s=t=0$ one gets (for any
$n$)
\begin{multline}\label{xx}
D_{E^{31}_{ii}E^{13}_{ii}}= \frac14\frac{\dd^2}{\dd
x^2_i}f(e^{a(x)})+\frac{\cosh 2x_i}{2\sinh 2x_i} \frac{\dd}{\dd
x_i}f(e^{ a(x)})\\
+\left(\frac{(\varkappa_{(2)}-\tilde{\varkappa}_{(1)})^2}{4\cosh^2
x_i}- \frac{(\tilde{\varkappa}_{(1)}-\varkappa_{(1)})^2}{4\sinh^2
x_i}-
\frac{(\varkappa_{(1)}-\varkappa_{(2)})}2(1+\frac{(\varkappa_{(1)}-\varkappa_{(2)})}2)
\right).
\end{multline}
In the algebra $\mathfrak k$ there is an identity
$[E_{ii}^{13},E_{ii}^{31}]= E^{11}_{ii}-E^{33}_{ii}$. Hence we
have $D_{E^{13}_{ii}E^{31}_{ii}}=
D_{E^{31}_{ii}E^{13}_{ii}}+\varkappa_{(1)}-\varkappa_{(2)}$.

The calculation of $D_{E^{23}_{ij}}D_{E^{32}_{ji}}f(e^{a(x)})$ in
the general  case is absolutely the same as in the case $n=1$,
$m=2$. For brevity we make this simplest calculation (in the
general case we only have to write indices $ij$ everywhere).

We need to translate the matrix $$ e^{a(x)}e^{sE^{23}}e^{tE^{32}}=
\begin{pmatrix}
\cosh x & t\sinh x &\sinh x\\
0 & 1+st & s\\
\sinh x & t\cosh x& \cosh x
\end{pmatrix},
$$
 by the left and right action of $K$ into the form $e^{a(y)}$, for some
$y\in\mathbb C$. That is we must find a representation of
$e^{a(x)}e^{sE^{23}}e^{tE^{32}}$ in the form
$e^{a(x)}e^{sE^{23}}e^{tE^{32}}=k_1 e^{a(y)} k_2$, $k_1,k_2\in K$,
$y\in\mathbb C$ . Applying a function $f\in F_{\vec{\varkappa}}$
to  both sides of this equation yields
\begin{multline}\label{rr}
f(e^{a(x)}e^{sE^{23}}e^{tE^{32}})=
\begin{pmatrix}
u^{\frac12}&
-\frac{t\tanh x}{\cosh x(1+st){\Delta}^2}&0\\
\frac{s}{{\Delta}\cosh x}&
u^{-1}&0\\
0&0&1
\end{pmatrix}f(e^{a(y)})\\
\times u^{\frac{\varkappa_{(1)}}2}
(1+st)^{\varkappa_{(1)}}\left(\frac{\cosh x}{\cosh y}
\right)^{\varkappa_{(1)}+\varkappa_{(2)}},
\end{multline}
where $u=1+\frac{st}{(1+st)\sinh^2 x}$, $\Delta=\sqrt{\tanh^2 x+\frac{st}
{(1+st)\cosh^2 x}}$, $y=\mathop{\rm arctanh}\nolimits \Delta$.

Taking the derivative $\frac{\dd^2}{\dd s\dd t}$ of (\ref{rr})
 at the point $s=t=0$ and using the
formulas for the right action of elements
$E^{11}=\tilde{\varkappa}_{(1)}$, $E^{22}=\varkappa_{(1)}$ one
gets the formula (already in the general case)
\begin{equation}\label{yy}
D_{E_{ij}^{32}E_{ji}^{23}}f(e^{a(x)})= \frac{\cosh x}{2\sinh
x}\frac{\dd}{\dd x_i}f(e^{a(x)}) +\left(
\frac{\varkappa_{(1)}-\varkappa_{(2)}}2
-\frac{\tilde{\varkappa}_{(1)}-\varkappa_{(2)}}{\sinh^2
x_i}\right) f(e^{a(x)}).
\end{equation}
Again we can calculate $D_{E_{ij}^{23}}D_{E_{ji}^{32}}$
through $D_{E_{ij}^{32}}D_{E_{ji}^{23}}$ by using the identity
 inside $\gotk$.

Using the equations $D_{E^{33}_{ii}}f=\varkappa_{(2)} f$,
$D_{E^{11}_{ii}}f=\varkappa_{(1)} f$, $D_{E^{33}_{ij}}f=
D_{E^{11}_{ij}}f=0$
 for $i\ne j$,
$f\in F_{\vec{\varkappa}}$ and  (\ref{DDij}), (\ref{xx}),
(\ref{yy}) results into the formula
\begin{multline*}
2R_{C_2}= \Delta_{A}+\sum_{i=1}^n \left(\frac{2(m-n)\cosh
x_i}{\sinh x_i}+ \frac{2\cosh 2x_i}{\sinh
2x_i}\right)\frac{\dd}{\dd x_i}
\\+
2\sum_{i=1}^n \sum_{j\ne i} \left(\frac{\cosh (x_i-x_j)}{\sinh
(x_i-x_j)}+ \frac{\cosh (x_i+x_j)}{\sinh
(x_i+x_j)}\right)\frac{\dd}{\dd x_i}\\
-2\varkappa_v(\varkappa_v+1)\sum_{i<j}\left(\frac{1}{\sinh^2
(x_i+x_j)}+ \frac{1}{\sinh^2 (x_i-x_j)}\right)\\
+\sum_{i=1}^n\frac{(\varkappa_{(2)}-\tilde{\varkappa}_{(1)})^2}{\cosh^2
x_i} -\frac{(\tilde{\varkappa}_{(1)}
-\varkappa_{(1)})(\tilde{\varkappa}_{(1)}-\varkappa_{(1)}+2(m-n))}{\sinh^2
x_i}
\\+
(\varkappa_{(1)}+\varkappa_{(2)})^2n+2(m-n)\varkappa_{(1)}^2.
\end{multline*}
Conjugating $R_{C_2}$ with $\delta_{\vec{s}}(x)$
and using a consequence of
the Weyl determinant formula for the $D_n$ root system:
\begin{multline*}
2\sum_{i<j}\left(\frac{1}{\sinh^2(x_i-x_j)}+
\frac{1}{\sinh^2(x_i+x_j)}\right)\\=
\sum_{i=1}^n\left(\sum_{j\ne i}\frac{\cosh(x_i-x_j)}{\sinh(x_i-x_j)}+
\frac{\cosh(x_i+x_j)}{\sinh(x_i+x_j)}\right)^2- (\rho_{0,1,0},
\rho_{0,1,0}),
\end{multline*}
 one gets the formula
(\ref{CMf}).
\end{proof}

\subsection{The branching rules}Let us recall the  branching rules  for
the inclusion $K\subset G$ (see \cite{King}). The construction is
based
 on the Littlewood-Richardson rule \cite{LR,McD} which
 deals with  partitions (and their diagrams).
A partition is a sequence of positive integer numbers
$\lambda\in\mathbb Z_{+}^r$ such that
$\lambda_1\ge\lambda_2\ge\dots\ge\lambda_r\ge 0$, $r$ is called
the length of the partition. Set
$|\lambda|=\sum_{i=1}^r\lambda_i$. The diagram of the partition is
the set of points $(i,j)\in \mathbb Z_{+}^n$ such that $1\le j\le
\lambda_i$. It is more convenient to replace the points by squares
(or boxes). We write $\mu\subset \lambda$ for the partitions if
and only if $\mu_i\le \lambda_i$ for all $i$. For a partition
$\lambda$ of the length less or equal $r$ let $L_\lambda$ be the
corresponding finite dimensional irreducible $GL(r)$
representation of the highest weight $\lambda$.

In these notations the branching rule has the form: $$
L_\lambda|_{K}=\sum_{\zeta}L_\zeta\boxtimes
\sum_{\tau}c^\lambda_{\zeta\tau} L_\tau, $$ where
$\lambda,\tau,\zeta$ are  partitions,
 and  $c^\lambda_{\zeta\tau}$ is a
non-negative integer coefficient given by the
Littlewood-Richardson rule. This coefficient is called the
Littlewood-Richardson number.

Let us recall the Littlewood-Richardson rule. For this purpose I
need
 some basic
combinatorial definitions. The set theoretic difference
$\theta=\lambda \setminus \mu$, $\mu\subset\lambda$ is called a
{\it skew diagram}, and $|\theta |=|\lambda|-|\mu|$. A skew
diagram is a {\it horizontal strip} if and only if it has at most
one square in each column. A {\it tableau}  $T$ is the sequence of
partitions (diagrams)
$\mu=\lambda^{(0)}\subset\lambda^{(1)}\subset\dots
\subset\lambda^{(r)}=\lambda$ such that each of the skew  diagrams
$\theta^{(i)}=\lambda^{(i)}-\lambda^{(i-1)}$ ($1\le i\le r$) is a
horizontal strip. Graphically, $T$ may be described by the
numbering each square of the skew diagram $\theta^{(i)}$
 with the number $i$. The numbers inserted in $\lambda-\mu$ must
increase strictly down each column and weakly from left to right
along each row. The skew diagram $\lambda-\mu$ is called  the
shape of tableau $T$, and the sequence
$(|\theta^{(1)}|,\dots,|\theta^{(r)}|)$ is called the weight of
$T$.

Let $T$ be a tableau. From $T$ one can derive a word $w(T)$ by
reading the symbols in $T$ from right to left in successive rows,
starting with the top row. A word $w=a_1a_2\dots a_N$ in the
symbols $1,2,\dots,n$ is said to be  a {\it lattice permutation}
if for $1 \le r\le N$ and $1\le i\le n-1$, the number of
occurrences of  $i$ in $a_1a_2\dots a_r$ is not less than the
number of the  occurrences of $i+1$.

\begin{picture}(55,35)
\put(0,0){\line(1,0){10}} \put(0,10){\line(1,0){10}}
\put(0,0){\line(0,1){10}} \put(10,0){\line(0,1){10}}
\put(3,1){$2$}

\put(10,0){\line(1,0){10}} \put(10,10){\line(1,0){10}}
\put(10,0){\line(0,1){10}} \put(20,0){\line(0,1){10}}
\put(13,1){$3$}

\put(10,10){\line(1,0){10}} \put(10,20){\line(1,0){10}}
\put(10,10){\line(0,1){10}} \put(20,10){\line(0,1){10}}
\put(13,11){$1$}

\put(20,10){\line(1,0){10}} \put(20,20){\line(1,0){10}}
\put(20,10){\line(0,1){10}} \put(30,10){\line(0,1){10}}
\put(23,11){$2$}

\put(20,20){\line(1,0){10}} \put(20,30){\line(1,0){10}}
\put(20,20){\line(0,1){10}} \put(30,20){\line(0,1){10}}
\put(23,21){$1$}

\put(30,20){\line(1,0){10}} \put(30,30){\line(1,0){10}}
\put(30,20){\line(0,1){10}} \put(40,20){\line(0,1){10}}
\put(33,21){$1$}

\put(40,20){\line(1,0){10}} \put(40,30){\line(1,0){10}}
\put(40,20){\line(0,1){10}} \put(50,20){\line(0,1){10}}
\put(43,21){$1$}

\put(60,20){For example the word $w(T)$ for the tableau from the
picture } \put(60,5){is $1112132$. It is an example of the lattice
permutation.}
\end{picture}

\begin{LRrule} Let $\lambda,\mu,\nu$ be partitions. Then
$c^{\lambda}_{\mu\nu}$
is zero unless
$\mu\subset\lambda$, $\nu\subset\lambda$, $|\mu|+|\nu|=|\lambda|$ and
for $\mu,\nu\subset\lambda$,  $|\mu|+|\nu|=|\lambda|$ it is
 equal to the number of tableaux $T$ of the shape $\lambda-\mu$ and weight
$\nu$ such that $w(T)$ is a lattice permutation.
\end{LRrule}
 The genuine definition of the Littlewood-Richardson number through Schur
functions \cite{McD} implies that  $c^{\lambda}_{\nu\mu}=c^{\lambda}_{\mu\nu}$.

Let $l$ be such a big integer  that $\varkappa_{(i)}+l>0$,
$\tilde{\varkappa}_{(1)}+l>0$ and the shifted highest weights
$\lambda'=\lambda+l1^{r}$ form a partition. Below we use the
superscript prime for the shifted objects.

\begin{rema}\label{l+} $c^{\lambda}_{\nu\mu}=c^{\lambda'}_{\nu'\mu'}$
\end{rema}

The last remark allows  us to define the coefficient
 $c^{\lambda}_{\nu\mu}$, when
one of the  $\lambda,\nu,\mu$ is not a partition, by the formula
from the remark. Below we suppose that all weight are shifted and
the superscript prime is suppressed.

\subsection{Combinatorial proofs}Now we use  the
Littlewood-Richardson rule to prove
lemmas~\ref{1dW},\ref{S},\ref{s1}.

\begin{proof}[Proof of  lemma~\ref{1dW}]
Every irreducible representation $L_{\mu}$ of $GL(m-1)$, where
$\mu\in P_+(GL(m-1))$ such that $\lambda_1\ge
\mu_1\ge\lambda_2\ge\mu_2\ge\lambda_3\ge\dots
\ge\mu_{m-1}\ge\lambda_{m}$ is contained at the restriction
$L_{\lambda}|_{GL(m-1)}$ exactly once ( see e.g.\cite{Zh} page
186). Applying this statement  $m-n$ times to $\lambda=
(\tilde{\varkappa}_{(1)}-\varkappa_{(1)})1^n+\varkappa_{(1)}1^m$
we get the first part of the lemma (e.g. $\dim
\tilde{W}_{\vec{\varkappa}}=1$).

Now let us use the Littlewood-Richardson rule for the restriction
from the group $K_{(1)}=GL(m)$ to the group $\tilde{K}_{(1)}\times
M= GL(n)\times GL(m-n)$. To prove the second part of the lemma one
must find all partitions $\nu$ of  length less or equal $n$ such
that $c^{\lambda}_{\nu,\varkappa_{(1)}1^{m-n}}\ne 0$. That is one
must find all fillings of $\lambda\setminus \varkappa_{(1)}
1^{m-n}$ by $1,\dots,n$ such that the result is the tableau
satisfying the lattice permutation condition. The first part of
the lemma says that if such filling exists then it is unique.

Let us construct this filling. Let us fill the last
$\tilde{\varkappa}_{(1)}-\varkappa_{(1)}$ squares of the $i$-th
row ($i=1,\dots,n$) by $i$ and the first $\varkappa_{(1)}$ squares
of the $(m-n+i)$-th row ($i=1,\dots,n$) by $i$. One can check that
the resulting tableau satisfies the lattice permutation condition
and has the weight $\nu= \tilde{\varkappa}_{(1)}1^n$.

\begin{picture}(35,35)
\put(0,0){\line(1,0){10}} \put(0,10){\line(1,0){10}}
\put(0,0){\line(0,1){10}} \put(10,0){\line(0,1){10}}
\put(3,1){$2$}

\put(0,10){\line(1,0){10}} \put(0,20){\line(1,0){10}}
\put(0,10){\line(0,1){10}} \put(10,10){\line(0,1){10}}
\put(3,11){$1$}

\put(10,10){\line(1,0){10}} \put(10,20){\line(1,0){10}}
\put(10,10){\line(0,1){10}} \put(20,10){\line(0,1){10}}
\put(13,11){$2$}

\put(20,10){\line(1,0){10}} \put(20,20){\line(1,0){10}}
\put(20,10){\line(0,1){10}} \put(30,10){\line(0,1){10}}
\put(23,11){$2$}

\put(10,20){\line(1,0){10}} \put(10,30){\line(1,0){10}}
\put(10,20){\line(0,1){10}} \put(20,20){\line(0,1){10}}
\put(13,21){$1$}

\put(20,20){\line(1,0){10}} \put(20,30){\line(1,0){10}}
\put(20,20){\line(0,1){10}} \put(30,20){\line(0,1){10}}
\put(23,21){$1$}

\put(40,20){For example on the picture we drew the filling
corresponding to}

\put(40,5){the case $m=3$, $n=2$, $\tilde{\varkappa}_{(1)}=3$,
$\varkappa_{(2)}=1$.}
\end{picture}

\end{proof}

To prove lemma~\ref{S} we must calculate
$c^{\lambda}_{\varkappa_{(1)} 1^m, \varkappa_{(2)}1^n}$. Let
$T=\{\varkappa_{(1)}1^m=\lambda^{(0)} \subset
\dots\subset\lambda^{(n)}=\lambda\}$ be a tableau contributing to
the Littlewood-Richardson number $c^{\lambda}_{\varkappa_{(1)}
1^m, \varkappa_{(2)}1^n}$. Then all the squares in the  $i$-th row
are labeled by $i$,  and the $(m+i)$-th row may contain only the
symbols $i,\dots,n$ ($1\le i\le n$). Also the horizontal strip
condition implies that $T$ can not contain any label at the
$(n+j)$-th row ($j=1,\dots,m-n$). Hence it implies
$\lambda_{n+j}=\varkappa_{(1)}$. We use following notations for
the labeling of the last $n$ rows: the number of occurrences of
the symbol $i$ at the $(m+j)$-th row is equal to $\mu^i_j$.

\begin{picture}(100,100)
\put(0,0){\line(1,0){10}} \put(0,10){\line(1,0){10}}
\put(0,0){\line(0,1){10}} \put(10,0){\line(0,1){10}}
\put(3,1){$4$}

\put(0,10){\line(1,0){10}} \put(0,20){\line(1,0){10}}
\put(0,10){\line(0,1){10}} \put(10,10){\line(0,1){10}}
\put(3,11){$3$}

\put(10,10){\line(1,0){10}} \put(10,20){\line(1,0){10}}
\put(10,10){\line(0,1){10}} \put(20,10){\line(0,1){10}}
\put(13,11){$4$}

\put(0,20){\line(1,0){10}} \put(0,30){\line(1,0){10}}
\put(0,20){\line(0,1){10}} \put(10,20){\line(0,1){10}}
\put(3,21){$2$}

\put(10,20){\line(1,0){10}} \put(10,30){\line(1,0){10}}
\put(10,20){\line(0,1){10}} \put(20,20){\line(0,1){10}}
\put(13,21){$3$}

\put(20,20){\line(1,0){10}} \put(20,30){\line(1,0){10}}
\put(20,20){\line(0,1){10}} \put(30,20){\line(0,1){10}}
\put(23,21){$4$}

\put(30,20){\line(1,0){10}} \put(30,30){\line(1,0){10}}
\put(30,20){\line(0,1){10}} \put(40,20){\line(0,1){10}}
\put(33,21){$4$}

\put(0,30){\line(1,0){10}} \put(0,40){\line(1,0){10}}
\put(0,30){\line(0,1){10}} \put(10,30){\line(0,1){10}}
\put(3,31){$1$}

\put(10,30){\line(1,0){10}} \put(10,40){\line(1,0){10}}
\put(10,30){\line(0,1){10}} \put(20,30){\line(0,1){10}}
\put(13,31){$2$}

\put(20,30){\line(1,0){10}} \put(20,40){\line(1,0){10}}
\put(20,30){\line(0,1){10}} \put(30,30){\line(0,1){10}}
\put(23,31){$3$}

\put(30,30){\line(1,0){10}} \put(30,40){\line(1,0){10}}
\put(30,30){\line(0,1){10}} \put(40,30){\line(0,1){10}}
\put(33,31){$3$}

\put(40,40){\line(1,0){10}} \put(40,50){\line(1,0){10}}
\put(40,40){\line(0,1){10}} \put(50,40){\line(0,1){10}}
\put(43,41){$4$}

\put(40,50){\line(1,0){10}} \put(40,60){\line(1,0){10}}
\put(40,50){\line(0,1){10}} \put(50,50){\line(0,1){10}}
\put(43,51){$3$}

\put(40,60){\line(1,0){10}} \put(40,70){\line(1,0){10}}
\put(40,60){\line(0,1){10}} \put(50,60){\line(0,1){10}}
\put(43,61){$2$}

\put(50,60){\line(1,0){10}} \put(50,70){\line(1,0){10}}
\put(50,60){\line(0,1){10}} \put(60,60){\line(0,1){10}}
\put(53,61){$2$}

\put(60,60){\line(1,0){10}} \put(60,70){\line(1,0){10}}
\put(60,60){\line(0,1){10}} \put(70,60){\line(0,1){10}}
\put(63,61){$2$}

\put(40,70){\line(1,0){10}} \put(40,80){\line(1,0){10}}
\put(40,70){\line(0,1){10}} \put(50,70){\line(0,1){10}}
\put(43,71){$1$}

\put(50,70){\line(1,0){10}} \put(50,80){\line(1,0){10}}
\put(50,70){\line(0,1){10}} \put(60,70){\line(0,1){10}}
\put(53,71){$1$}

\put(60,70){\line(1,0){10}} \put(60,80){\line(1,0){10}}
\put(60,70){\line(0,1){10}} \put(70,70){\line(0,1){10}}
\put(63,71){$1$}

\put(70,70){\line(1,0){10}} \put(70,80){\line(1,0){10}}
\put(70,70){\line(0,1){10}} \put(80,70){\line(0,1){10}}
\put(73,71){$1$}

\put(90,80){ The drawn tableau $T$ corresponds to the case $n=4$,}

\put(90,65){ $m=4$, $\varkappa_{(1)}=4$, $\varkappa_{(2)}=1$. The
diagram $\lambda$ has the }

\put(90,50){shape: $\lambda=(8,7,5,5,4,4,2,1)$. The labeling of
the}

\put(90,35){last $n$ rows encoded by $\mu$: $\mu^1_1=1$,
$\mu^2_1=1$, $\mu^3_1=2$,}

\put(90,15){ $\mu^4_1=0$; $\mu^2_2=1$, $\mu^3_2=1$, $\mu^4_2=2$;
$\mu_3^3=1$, $\mu_3^4=1$;}

\put(90,0){ $\mu^4_4=1$.   The word $w(T)$ is a lattice
permutation. }

\end{picture}

\begin{lem}\label{skl}
If $c^{\lambda}_{\varkappa_{(1)}1^m,\varkappa_{(2)}1^n}\ne 0$ then
\begin{equation}\label{ll}
\lambda_{j}+\lambda_{m+n+1-j}=\varkappa_{(1)}+ \varkappa_{(2)},
\qquad (j=1,\dots,n)
\end{equation}
\begin{equation}\label{m-n}
\lambda_{n+i}=\varkappa_{(1)} \qquad (j=1,\dots,m-n)
\end{equation}
\begin{equation}
\label{ln} \lambda_n\ge\varkappa_{(1)}.
\end{equation}
 In this case
$c^{\lambda}_{\varkappa_{(1)} 1^m,\varkappa_{(2)}1^n}=1$ and the
labeling of the corresponding tableau is given by formula:
\begin{equation}\label{mu}
\mu^i_j=\lambda_{i-j}-\lambda_{i-j+1},
\end{equation}
where $i\ge j$ and $\lambda_0=\varkappa_{(1)}+\varkappa_{(2)}$.
\end{lem}
\begin{proof}
 Let us remark that the reasoning before the lemma proves formula
(\ref{m-n}).
 Hence the $n+j$-th row $(j=1,\dots,m-n)$ plays no role and we consider only
the case $n=m$ in the proof.

We prove the claim by the induction in $n$.

For $n=1$ the claim is obvious.

Let $T=\{\lambda^{(0)}\subset\dots\subset\lambda^{(s+1)}\}$ be  a
tableau contributing to $c^{\lambda^{(s+1)}}_{\varkappa_{(1)}
1^{s+1},\varkappa_{(2)}1^{s+1}}$
 and $\tilde T$ is the tableau
obtained from $T$ by deleting the boxes with $s+1$.
 Tableau $\tilde T$ contains
only the symbols $1,\dots,s$ and is of the shape
$\tilde{\lambda}\setminus\varkappa_{(1)} 1^s$ for partition
$\tilde{\lambda}$, and $\tilde{\lambda_i}=\lambda_i$
($i=1,\dots,s$). Furthermore, $\tilde{T}$ contributes to the
Littlewood-Richardson number $c^{\tilde{\lambda}}_{\varkappa_{(1)}
1^s,\varkappa_{(2)}1^s}$.  Hence for $\tilde T$
(\ref{ll})-(\ref{mu}) hold by the induction hypothesis. That is,
the numbers $\mu^i_j$, $s\ge i\ge j$, for given $\lambda$ are
uniquely determined by (\ref{mu}) and we only need to find
$\mu^{s+1}_j$, $j=1,\dots,s+1$. The horizontal strip condition and
the induction hypothesis imply
\begin{equation*}
\sum_{i=1}^{s+1}\mu_j^i\le\sum^s_{i=1}\mu_{j-1}^i=\lambda_0-\lambda_{s+2-j},
\end{equation*}
for $j=2,\dots,s+1$.  The induction hypothesis for
$j=1,\dots,s+1$, implies
\begin{multline*}
\sum_{i=1}^{s+1}\mu_j^i=\tilde{\lambda}_{s+j}+\mu^{s+1}_j=
\varkappa_{(1)}+\varkappa_{(2)}-\tilde{\lambda}_{s+1-j}+\mu_j^{s+1}\\=
\varkappa_{(1)}+\varkappa_{(2)}-\lambda_{s+1-j}+\mu_j^{s+1},
\end{multline*}
where $\tilde{\lambda}_{2s}=0$. Hence for $j=2,\dots,s+1$ we have
\begin{equation}\label{hs}
\lambda_{s+2-j}\le \lambda_{s+1-j}-\mu_{j}^{s+1}.
\end{equation}
The lattice permutation condition for $T$ implies
\begin{equation}\label{lp}
\lambda_{s+1}+\mu_1^{s+1}\ge \lambda_s.
\end{equation}
Adding these $s+1$ inequalities,  one gets
$\sum_{j=1}^{s+1}\lambda_j
+\sum_{j=1}^{s+1}\mu_j^{s+1}\le\sum^s_{j=0}\lambda_j$, but the
weight condition for $T$ implies that the last inequality is an
equality. Hence (\ref{hs}), (\ref{lp}) are also equalities and
they imply (\ref{mu}). Thus we proved that
$c^{\lambda}_{\varkappa_{(1)}1^{s+1},\varkappa_{(2)}1^{s+1}}\ne 0$
implies (\ref{mu}) for $n=s+1$. One can easily check that $T$
defined by (\ref{mu}) contributes to the Littlewood-Richardson
number
$c^{\lambda}_{\varkappa_{(1)}1^{s+1},\varkappa_{(2)}1^{s+1}}$.
\end{proof}
 Lemma~\ref{skl} is equivalent to lemma~\ref{S}.

Let us denote by the symbols $v$ and $w$
 the partitions $v=
(\varkappa_{(1)}+\varkappa_{(2)}-\tilde{\varkappa}_{(1)})1^n+\varkappa_vne_1$,
where $e_1=(1,0,\dots,0)$, and
$w=\tilde{\varkappa}_{(1)}1^m+\varkappa_{(1)}1^n$.

\begin{lem}\label{sl}
 If $\lambda$ is a partition of the length less or equal $n+m$
 such that equalities (\ref{ll}) and
inequalities
 \begin{equation}\label{ineq1}
\lambda_n\ge\tilde{\varkappa}_{(1)},
\end{equation}
\begin{equation}\label{ineq2}
\lambda_i-\lambda_{i+1}\ge \varkappa_v,
\end{equation}
where $i=1,\dots,n-1$,
  hold. Then
$c^{\lambda}_{w,v}=1$ and the labeling of the corresponding
tableau is given by the formulas:
\begin{equation}\label{eq1}
\mu_j^i=\lambda_{i-j}-\lambda_{i-j+1}, \mbox{  for } i\ge j>1,
\end{equation}
\begin{equation}\label{eq2}
\mu^i_1=\lambda_{i-1}-\lambda_i-\varkappa_v, \mbox{  for } i>1,
\end{equation}
\begin{equation}\label{eq3}
\mu^1_1=(n-1)\varkappa_v+\varkappa_{(1)} -\lambda_1,
\end{equation}
where $\lambda_0=\varkappa_{(1)}+\varkappa_{(2)}$ and  $\mu^i_j=0$
for $i<j$.
\end{lem}
\begin{proof} Let $\lambda$ be a partition of the length less or equal $m+n$,
satisfying  (\ref{eq1}), (\ref{eq2}). Let $T$ be a tableau of the
shape $\lambda\setminus w$ such that the $i$-th row is filled by
the symbol $i$ ($i=1,\dots,n$) and the $(m+j)$-th string contains
$\mu^i_j$ symbols $i$, where $\mu_j^i$ are given by formulas
(\ref{eq1})-(\ref{eq3}). One can check that this tableau
contributes to $c^{\lambda}_{w, v}$.

Now we will prove that if $\lambda$ satisfies the conditions from
the lemma then
  $c^{\lambda}_{w, v}=1$, and equations (\ref{eq1})-(\ref{eq3}) hold.
  We will do it
by the induction. There is no difference in reasonings in the case
$m=n$ and in the case $m>n$, and we consider only the first case.
In this case $w=\tilde{\varkappa}_{(1)}1^n$.

For $n=2$ the claim is obvious.

Now let $T=\{\lambda^{(0)}\subset\dots\subset\lambda^{(s+1)}=\lambda\}$
be a skew tableau
satisfying inequalities (\ref{ineq1})-(\ref{ineq2}) for $n=s+1$
 which contributes to  $c^{\lambda}_{(\tilde{\varkappa}_{(1)})1^{s+1},
 v'}$, where
$v'=s\varkappa_v e_1+(\varkappa_{(1)}+\varkappa_{(2)}-
\tilde{\varkappa}_{(1)})1^{s+1}$. Let us remove $\varkappa_v$
symbols $1$ from the $(s+2)$-th row of $T$, delete the $(s+1)$-th
and $2(s+1)$-th rows and all boxes with $s+1$ from $T$. Then one
gets a tableau $\tilde T$ of some skew shape
$\tilde{\lambda}\setminus (\tilde{\varkappa}_{(1)})1^s$, which
contributes to the Littlewood-Richardson number
$c^{\tilde{\lambda}}_{(\tilde{\varkappa}_{(1)})1^s, v''}$,
 $v''=\varkappa_v(s-1)e_1+(\varkappa_{(1)}+\varkappa_{(2)}-
\tilde{\varkappa}_{(1)}) 1^s$. Moreover, for $\tilde T$ the
inequalities (\ref{ineq1})-(\ref{ineq2}) hold, hence by the
induction hypothesis equations (\ref{eq1})-(\ref{eq3}) hold for
$\tilde T$. Thus we found  $\mu_i^j$, $i\le s$ for $T$, and we
only need to find $\mu^{s+1}_j$. But we know $\lambda_i$,
$i=1,\dots,2(s+1)$ and
 $\mu_i^j$, $i\le s$, hence we can calculate $\mu^{s+1}_j$.
\end{proof}
  Remark~\ref{l+} and lemma~\ref{sl}
imply lemma~\ref{s1}.

\subsection{Proof of lemma~\ref{Wact} and the asymptotic estimate}
\begin{proof}[Proof of  lemma~\ref{Wact}]
Let $w\in W$ to be an element of the Weyl group and $\theta_w\in
\tilde{K}_{(1)}\times K_{(2)}$ such that  $\theta_w^{11}=t(w)$,
$\theta_w^{33}=\hat{t}(q(w))$, here $\hat{t}$ is the standard
embedding $S_n\hookrightarrow GL(n,\mathbb Z)$.
 Then for $f\in  F_{\vec{\varkappa}}$:
\begin{multline*} f(e^{a(w(x))})=f(\theta_we^{a(x)}\theta_{w^{-1}})
=\theta^wf(e^{a(x)})\chi_0^{\varkappa_{(1)}}\chi_1^{\varkappa_{(2)}}(w)\\
=\chi_0^{\tilde{\varkappa}_{(1)}}\chi_1^{\varkappa_{(2)}+\varkappa_{(1)}
-\tilde{\varkappa}_{(1)}-\varkappa_v}(w)
f(e^{a(x)})\chi_0^{\varkappa_{(1)}}\chi_1^{\varkappa_{(2)}}(w)\\=
\chi_0^{\kappa_1+\kappa_3}\chi_1^{\kappa_1+\kappa_2+\kappa_3}(w)f(e^{a(x)}),
\end{multline*}
here the third equality follows from  lemma~\ref{1dW} and the fact
that $q(w)$ acts on $\tilde{V}_{\vec{\varkappa}}$ by
$\det(\hat{t}(q(w)))^{\varkappa_{(2)}+\varkappa_{(1)}
-\tilde{\varkappa}_{(1)}-\varkappa_v}$.

 The element $e^{a(\pi i e_j)}$ belongs to the subgroup $K$,
 hence $f(e^{a(x+\pi i e_j)})=
f(e^{a(x)}e^{a(\pi i
e_j)})=(-1)^{\varkappa_{(1)}+\varkappa_{(2)}}f(e^{a(x)})
=(-1)^{\kappa_1}f(e^{a(x)})$.
\end{proof}

We can estimate the asymptotic behavior at the infinity of a
matrix element of the representation $L_\lambda$. One says that an
asymptotic estimate $f(x)\lesssim g(x)$ holds in the sector,
$x_1>x_2>\dots>x_n$ if  for any $y$ from the sector the limit
$\lim_{t\to +\infty}\frac{f(ty)}{g(ty)}$ is finite.

\begin{lem}\label{as}
Let $v\in L_\lambda$, $u\in L^*_\lambda$, then at the sector
$x_1>x_2>\dots>x_n$ we have an asymptotic estimate
\begin{equation}\label{asy}
\langle u ,e^{h(x)}v\rangle\lesssim
\prod_{i=1}e^{x_i(\lambda_i-\lambda_{m+n+1-i})}.
\end{equation} Moreover
\begin{equation}\label{lim}
\lim_{t\to +\infty}
\langle u ,e^{a(x)}v\rangle
\prod_{i=1}e^{-x_i(\lambda_i-\lambda_{m+n+1-i})}=
\langle u,v_\lambda\rangle\langle v,v_\lambda^*\rangle,
\end{equation}
where $v_\lambda$ and $v_\lambda^*$ are the highest and lowest
weight vectors of $L_\lambda$ and $L^*_\lambda$, respectively, and
$\langle v_\lambda,v_\lambda^* \rangle=1$.
\end{lem}
\begin{proof} We use notations from the subsection~\ref{restrr}.
There is a Cartan subgroup  $H\subset G$, $A=\exp(\mathfrak p)\cap
H$ and $e^{a(x)}=h(x,0,-x)$.
 The highest weight of $L_\lambda$
with respect to $ H$ is equal to $\lambda$, and all other
extremal  weights
are of the form $w(\lambda)$, $w\in W$. Obviously, for  proving
the claim it is enough to prove the estimate for the case when
 $u,v$ are extremal weight vectors. If $u,v$ are the extremal weight vectors then
   $\langle u, e^{a(x)} v\rangle\sim
 \prod_{i=1}e^{(w(\lambda)_i-w(\lambda)_{m+n+1-i})x_i}$, for $x\sim
 \infty$, and  obviously  in the sector $x_1>\dots>x_n$  the asymptotic
 estimate (\ref{asy}) holds.
\end{proof}

\subsection{Proof of lemma~\ref{dd}}
For $f\in F_{\vec{\varkappa}}$ the equality $f(e^{a(x)})=f(s(x))$
holds, where $s(x)\in G$ such that $s(x)^{11}=s(x)^{33}=1$,
$s(x)^{13}=s(x)^{31}=diag(z_1,\dots,z_n)$, $s(x)^{22}=1$,
$s(x)^{23}=s(x)^{32}=s(x)^{21}=s(x)^{12}=0$,
 $z_i=\tanh x_i$,
$i=1,\dots,n$.

\begin{lem}\label{diveij} For any $1\le i<j\le n$ the function
$(E^{11}_{ij}\pm
E^{22}_{ij})^m\frac{f(s(x))}{(z_i\pm z_j)^m}$ is regular at the generic
point $x$ such that $\sinh(x_i\pm x_j)=0$.
\end{lem}
\begin{proof}
For $f\in F_{\vec{\kappa}}$ the following equation holds $$
e^{y(E^{11}_{ij}+E^{33}_{ij})}f(s(x))= f( e^{y(E_{ij}^{11}+
E^{33}_{ij})}s(x)e^{-y(E_{ij}^{11}+ E^{33}_{ij})})=f(m(x,y)),$$
where $m^{11}(x,y)=m^{33}(x,y)=1$, $y\in\mathbb C$ and
\begin{multline*}
m^{13}=m^{31}=e^{yE_{ij}}ze^{-yE_{ij}}=(1+yE_{ij})
z
(1-yE_{ij})\\=z+y[E_{ij},z]=
z+y(z_i-z_j)E_{ij},
\end{multline*}
where $z=diag(z_1,\dots,z_n)$.
 Hence the function
$e^{\frac{t(E^{11}_{ij}+E^{33}_{ij})}{z_i-z_j}}f(s(x))$ is regular
at the generic point of $z_i=z_j$. By taking the derivative
$\frac{d}{dt}$ at the point $t=0$ one gets  the claim for $m=1$.
Iterating this procedure we obtain the proof in the case of the
minus sign. In case of the plus sign one  can proceed analogously
by considering $e^{y(E_{ij}^{11}-E_{ij}^{33})}f(s(x))=
f(e^{y(E_{ij}^{11}- E^{33}_{ij})}s(x)e^{-y(E_{ij}^{11}-
E^{33}_{ij})}).$
\end{proof}

\begin{proof}[Proof of lemma~\ref{dd}]
It is easy to see that all the weight subspaces (with respect to
action of $T$) of $\tilde{W}_{\vec{\varkappa}}\boxtimes
V_{\vec{\varkappa}}$ are one-dimensional. Hence for $v\in
\tilde{W}_{\vec{\varkappa}}\boxtimes \tilde{V}{\vec{\varkappa}}$,
  we
have $(E^{11}_{ij}\pm E^{33}_{ij})^lv=(\pm E^{33}_{ij})^lv\ne 0$,
for $l\le \kappa_2$, $i\ne j$ and zero otherwise. Using the last
remark and the simple trigonometric identity $z_i\pm
z_j=\frac{\sinh(x_i\pm x_j)}{\cosh x_i\cosh x_j}$, one derives
from  lemma~\ref{diveij} the divisibility of $f(e^{a(x)})$ by
$\delta_{0,\kappa_2,0}$.

In the rest part of the proof we consider only the case $m=n$
because the case $m>n$ is absolutely analogous.
 For $f\in F_{\vec{\varkappa}}$ following equation holds
\begin{multline*}
\delta^{-1}_{\kappa_1,0,\kappa_3}(x)f(e^{a(x)})=\begin{pmatrix}\cosh^{-1}
x & 0\\0&\sinh x
\end{pmatrix}
f(e^{a(x)})
\begin{pmatrix}1
& 0
\\0&\tanh^{-1} x
\end{pmatrix}\\=
f\left(\begin{pmatrix}\cosh^{-1} x & 0\\0&\sinh x
\end{pmatrix}
e^{a(x)}
\begin{pmatrix}1
& 0\\0&\tanh^{-1} x\end{pmatrix}\right)= f\left(\begin{pmatrix}1 &
 1\\\sinh^2 x&\cosh^2 x
\end{pmatrix}\right)
\end{multline*}
The right hand side of the last equation is  regular at
the generic point $x$ such
that $\delta_{0,0,1}(x)=0$, hence the left hand side is also regular. Thus
$f(e^{a(x)})$ is divisible by $\delta_{\kappa_1,0,\kappa_3}$.
\end{proof}

\subsection{Proofs of the corollaries}
\begin{proof}[Proof of corollary~\ref{degr}]
If $\lambda=\tau(\mu+ \rho_{\vec{\kappa}})$ then
$\lambda_i-\lambda_{m+n+1-i}=2(\mu+\rho_{\vec{\kappa}})_i$. In the
sector $x_1>\dots>x_n$ $\delta_{\vec{\kappa}}$, has an asymptotic
behavior $\delta_{\vec{\kappa}}(x)\sim
e^{2(x,\rho_{\vec{\kappa}})}$. Hence in this sector the asymptotic
estimate $\Psi_\mu(x)/\delta_{\vec{\kappa}}(x)\lesssim
e^{2(x,\mu)}$ holds. Together with lemma~\ref{dd} it gives the
proof of the first item of the
 corollary.

The second item immediately follows from the first one. Indeed, if for some
$\mu\in\mathbb Z^n$, $\mu_1\le\mu_2\le\dots\le \mu_n$,
 $\lambda=\tau(\mu+
\rho_{\vec{\kappa}})$ belongs  to $P^{\vec{\varkappa}}$ then
$\Psi_\mu\ne 0$. Hence, by the previous item $\mu\ge 0$.
\end{proof}

\begin{proof}[Proof of corollary~\ref{dne0}]
Suppose that $d_{\mu\mu}=0$. Then there is an expansion
$\Psi_\mu/\delta_{\vec{\kappa}}=\sum_{\nu<\mu}c_{\mu\nu}
J^{\vec{r}}_{\nu}$. By the definition of the Jack polynomials
\begin{multline*}
(\Psi_\mu/\delta_{\vec{\kappa}},\tilde{L}_{\vec{r}}
\Psi_\mu/\delta_{\vec{\kappa}})_{\vec{r}}=
(\sum_{\nu<\mu}c_{\mu\nu} J^{\vec{r}}_{\nu}, \sum_{\nu<\mu}
(\nu+\rho_{\vec{r}},\nu+\rho_{\vec{r}})c_{\mu\nu}
J^{\vec{r}}_{\nu})_{\vec{r}}\\=
\sum_{\nu<\mu}(\nu+\rho_{\vec{r}},\nu+\rho_{\vec{r}}) c^2_{\mu\nu}
(J^{\vec{r}}_{\nu},J^{\vec{r}}_{\nu})_{\vec{r}}<
(\mu+\rho_{\vec{r}},\mu+\rho_{\vec{r}})
(\Psi_\mu/\delta_{\vec{\kappa}},\Psi_\mu/\delta_{\vec{\kappa}})_{\vec{r}}.
\end{multline*}
  Theorem~\ref{CM} and formula (\ref{c2}) yields
$\tilde{L}^{\vec{r}} \Psi_\mu/\delta_{\vec{r}}=
(\mu+\rho_{\vec{r}},\mu+\rho_{\vec{r}})\Psi_\mu/\delta_{\vec{r}}$.
\end{proof}

\subsection{Proof of the main theorem}
\begin{proof}[Proof of theorem~\ref{main}]
Corollary~\ref{dne0}, theorem~\ref{CM} and formula (\ref{c2})
imply the first item of the theorem.

The last items follow from the fact that  a $W$-invariant differential
operator is uniquely determined by its action  on the
space $\mathbb C[P]^W$ of $W$-invariant polynomials (see for
example page~16 of \cite{HS}).
Indeed,  formula (\ref{Cch})
 implies that the highest
term of  $R_{C_{2r}}$ has the form described at the theorem. As
$C_{2r}$ are pairwise commutative, hence $R_{C_{2r}}$ are also pairwise
commutative.

\end{proof}

\end{document}